\documentclass[a4paper,leqno,draft]{article}
\usepackage{amssymb,amsmath,amsfonts,amsthm,
mathrsfs}

\renewcommand{\theequation}{\thesection.\arabic{equation}}
\renewcommand{\theequation}{\thesection.\arabic{equation}}
\newtheorem{theorem}{Theorem}[section]
\newtheorem{corollary}[theorem]{Corollary}
\newtheorem{lemma}[theorem]{Lemma}
\newtheorem{proposition}[theorem]{Proposition}
\newtheorem{definition}[theorem]{Definition}
\theoremstyle{definition}
\newtheorem{remark}[theorem]{Remark}
\newtheorem{example}[theorem]{Example}

\newcommand{\wt}[1]{\widetilde{#1}}

\renewcommand{\S}{\mathscr{S}}

\newcommand{\mb}[1]{\ensuremath{\mathbb{#1}}}
\newcommand{\N}{\mb{N}}

\newcommand{\R}{\mb{R}}
\newcommand{\C}{\mb{C}}

\newcommand{\LL}{\mathcal{L}}

\newcommand{\G}{\ensuremath{{\cal G}}}

\newcommand{\Gc}{\ensuremath{{\cal G}_\mathrm{c}}}

\newcommand{\GS}{\G_{{\, }\atop{\hskip-4pt\scriptstyle\S}}\!}
\newcommand{\EM}{\ensuremath{{\cal E}_{M}}}

\newcommand{\Neg}{\mathcal{N}}

\newcommand{\Ginf}{\ensuremath{\G^\infty}}

\newcommand{\singsupp}{\mathrm{sing\, supp}}

\newfont{\bigmath}{cmr12 at 13pt}

\newfont{\grecomath}{cmmi12 at 15pt}

\newcommand{\val}{\mathrm{v}} 
\newcommand{\esp}{\mathrm{e}}

\newfont{\bl}{msbm10 scaled \magstep2}

\newcommand{\beq}{\begin{equation}}
\newcommand{\eeq}{\end{equation}}

\newcommand{\F}{\ensuremath{{\cal F}}}

\newcommand{\notmid}{\mid\kern-0.5em\not\kern0.5em}

\newcommand{\eps}{\varepsilon}

\newcommand{\Om}{\Omega}

\renewcommand{\Re}{\ensuremath{\text{Re}}}
\renewcommand{\Im}{\ensuremath{\text{Im}}}

\renewcommand{\Re}{\operatorname{Re}}
\renewcommand{\Im}{\operatorname{Im}}

\newcommand{\mB}{\mathcal{B}}
\newcommand{\mC}{\mathcal{C}}
\newcommand{\mL}{\mathcal{L}}
\newcommand{\mP}{\mathcal{P}}
\newcommand{\mQ}{\mathcal{Q}}
\newcommand{\mU}{\mathcal{U}}
\newcommand{\mF}{\mathcal{F}}
\newcommand{\mH}{\mathcal{H}}

\newcommand{\mS}{\mathcal{S}}
\newcommand{\mV}{\mathcal{V}}
\newcommand{\mW}{\mathcal{W}}

\begin{document}
\title{{\bf Closed graph and open mapping theorems for topological $\wt{\C}$-modules and applications}}
 \author{Claudia Garetto\footnote{Supported by FWF (Austria), grants T305-N13 and  Y237-N13, and TWF (Tyrol), grant
UNI-0404/305.} \\
Institut f\"ur Technische Mathematik,\\ 
Universit\"at Innsbruck, Austria\\
\texttt{claudia@mat1.uibk.ac.at}\\
}
\maketitle
\date{}
\maketitle

\begin{abstract}
We present closed graph and open mapping theorems for $\wt{\C}$-linear maps acting between suitable classes of topological and locally convex topological $\wt{\C}$-modules. This is done by adaptation of De Wilde's theory of webbed spaces and Adasch's theory of barrelled spaces to the context of locally convex and topological $\wt{\C}$-modules respectively. We give applications of the previous theorems to Colombeau theory as well to the theory of Banach $\wt{\C}$-modules. In particular we obtain a necessary condition for $\Ginf$-hypoellipticity on the symbol of a partial differential operator with generalized constant coefficients.
\end{abstract}

\setcounter{section}{-1}
\section{Introduction}
In the recent past the topological investigation of Colombeau algebras of generalized functions has motivated the development of a more general theory of locally convex and topological $\wt{\C}$-modules. This is based on new concepts as $\wt{\C}$-absorbent, balanced and convex subsets of a $\wt{\C}$-module $\G$ and makes use of a suitable $\wt{\C}$-adaptation of the classical notion of seminorm called ultra-pseudo-seminorm. The theory of locally convex and topological $\wt{\C}$-modules has been elaborated in \cite{Garetto:05a, Garetto:05b} by following, as a blueprint, the classical results on topological vector spaces concerning continuity of linear maps, boundedness, completeness, inductive and projective limit topologies etc. Particular attention has been given to the dual of a topological $\wt{\C}$-module $\G$, i.e. the set $\LL(\G,\wt{\C})$ of all continuous $\wt{\C}$-linear functionals on $\G$, and to the different ways of endowing $\LL(\G,\wt{\C})$ with a $\wt{\C}$-linear topology. 

The purpose of this paper is to provide closed graph and open mapping theorems for $\wt{\C}$-linear maps acting on suitable classes of locally convex and topological $\wt{\C}$-modules. This is obtained by working out a $\wt{\C}$-linear version of the classical theories of De Wilde's webs and Adasch's strings to employ on locally convex and topological $\wt{\C}$-modules respectively.  Since the most common Colombeau algebras are examples of locally convex topological $\wt{\C}$-modules these abstract theorems of functional analysis can now be applied to concrete problems in the Colombeau context and provide results of continuity and regularity. As a first application of the closed graph theorem we prove a necessary condition of $\Ginf$-hypoellipticity for partial differential operators with constant Colombeau coefficients and continuous dependence on the data for the solution of a Cauchy problem uniquely solvable in the Colombeau algebra $\G(\R^n)$.

We now describe the contents of the sections in more detail.

Section \ref{section_wilde} lays down an adaptation of De Wilde's theory of webbed spaces to the framework of locally convex topological $\wt{\C}$-modules. By replacing the classical real scalars with a special family of generalized numbers, as already done in many situations in \cite{Garetto:05a,Garetto:05b}, we introduce the notion of web of type $\mathcal{C}$ on a locally convex topological $\wt{\C}$-module $\G$. A first version of the closed graph theorem is obtained in Theorem \ref{theorem_closed} for a $\wt{\C}$-linear map acting from a Fr\'echet $\wt{\C}$-module $\G$ to a locally convex topological $\wt{\C}$-module $\mF$ which has an absolutely convex web of type $\mathcal{C}$. The validity of this result is extended by assuming that $\G$ is a Baire $\wt{\C}$-module in Theorem \ref{theorem_closed_Baire} and a Hausdorff ultrabornological $\wt{\C}$-module in Proposition \ref{prop_cor_ultrabor}. Some hereditary properties of webbed $\wt{\C}$-modules are examined in Subsection \ref{subsection_webs}. In particular we prove that the strict inductive limit of a sequence of locally convex topological $\wt{\C}$-modules equipped with a web of type $\mathcal{C}$ as well as its topological dual has a web of the same type.

Section \ref{section_open} is devoted to a first formulation of the open mapping theorem via application of the closed graph theorems provided in Section \ref{section_wilde}. In this way we obtain that a sequentially closed $\wt{\C}$-linear map from a locally convex topological $\wt{\C}$-module with a web of type $\mathcal{C}$ onto a Hausdorff ultrabornological $\wt{\C}$-module is open. By following \cite[Chapter 3, Section 17]{Horvath:66} we give a first example of an open mapping theorem which holds for $\wt{\C}$-linear maps acting between topological $\wt{\C}$-modules which are not necessarily locally convex but metrizable. The investigation of classes of topological $\wt{\C}$-modules without convexity assumptions by means of which to formulate closed graph and open mapping theorems is the main topic of Section \ref{section_adasch}.

Section \ref{section_adasch} elaborates an adaptation of Adasch's theory of strings for topological vector spaces to the setting of topological $\wt{\C}$-modules. In Subsection \ref{subsec_strings} after the basic definitions we deal with different $\wt{\C}$-linear topologies generated by directed families of strings and we make a comparison between them. The notion of barelledness for topological $\wt{\C}$-modules as well as the strong and the associated barrelled topologies on a $\wt{\C}$-module $\G$ is introduced in Subsection \ref{subsec_barrelled}. The core of Subsection \ref{subsec_closed_bar} is a version of the closed graph theorem which says that every closed $\wt{\C}$-linear map from a barrelled $\wt{\C}$-module $\G$ to a Fr\'echet $\wt{\C}$-module $\mF$ is continuous. In order to extend this theorem to a larger class of spaces $\mF$ we define the concept of infra-s $\wt{\C}$-module. This turns out to be the largest possible class of $\wt{\C}$-modules $\mF$ for which the closed graph theorem above is valid. Indeed, an infra-s $\wt{\C}$-module $\mF$ can be characterized by the property that every closed $\wt{\C}$-linear map from a barrelled $\wt{\C}$-module to $\mF$ is continuous. The final Subsection \ref{subsec_open_bar}  of Section \ref{section_adasch} deals with open mapping theorems and barrelled topological $\wt{\C}$-modules. We begin in Theorem \ref{theorem_open_mapping_bar} by proving that every surjective and continuous $\wt{\C}$-linear map from a Fr\'echet $\wt{\C}$-module to a Hausdorff barrelled topological $\wt{\C}$-module is open. By means of technical notions such as weakly singular $\wt{\C}$-linear maps and regular contractions we work out an extension of the previous theorem which allows to substitute the Fr\'echet $\wt{\C}$-modules with the larger class of s-$\wt{\C}$-modules. As for infra-s $\wt{\C}$-modules we prove that this choice is optimal in the sense that a Hausdorff topological $\wt{\C}$-module $\G$ is an s-$\wt{\C}$-module if and only if every weakly singular $\wt{\C}$-linear map from $\G$ onto a barrelled topological $\wt{\C}$-module is open. 

Some applications of the closed graph and the open mapping theorems to Colombeau theory as well as to the theory of Banach $\wt{\C}$-modules are collected in Section \ref{sec_application}. The recent investigation of the $\Ginf$-regularity properties of generalized differential and pseudodifferential operators in the Colombeau context \cite{Garetto:04, GGO:03, GH:05, GH:05b, HO:03, HOP:05} has provided several sufficient conditions of $\Ginf$-hypoellipticity, i.e. hypotheses on the generalized symbol of the operator $P(x,D)$ which allow to conclude that $u\in\Ginf(\Om)$ when $Pu\in\Ginf(\Om)$. So far, the research of necessary condition for $\Ginf$-hypoellipticity has been a long-standing open problem. 

The closed graph theorem stated for Fr\'echet $\wt{\C}$-modules as a particular case of Theorem \ref{theorem_closed_1} enables us to find for the first time a necessary condition for $\Ginf$-hypoellipticity on the symbol of a partial differential operator with generalized constant coefficients. The achievement of a result of this kind via functional analytic methods is a novelty within Colombeau theory and we strongly believe that closed graph theorem arguments will provide necessary conditions for other types of regularity such as Colombeau-analytic or Colombeau-Gevrey regularity.  A further immediate application of the closed graph theorem shows that when the Cauchy problem 
\[
u^{(n)}(t)=\sum_{i=1}^{n-1}a_i(t)u^{(i)}(t)+b(t),\qquad\quad u({t_0})=u_0,\ u'(t_0)=u_1,\,...,\,u^{(n-1)}(t_0)=u_{n-1}
\]
is uniquely solvable in $\G(\R)$ then the solution $u$ depends continuously on the inhomogeneity $b$ and the initial values $u_0, u_1,...,u_{n-1}$. Analogous continuous dependence on the data is obtained for the Colombeau solution of the linear wave equation problem.

The second part of Section \ref{sec_application} presents some theoretical applications of the closed graph theorem in the form of the isomorphism theorem to the abstract theory of Banach $\wt{\C}$-modules. In detail, we prove that in a Banach $\wt{\C}$-module two compatible ultra-pseudo-norms are equivalent and that a pointwise bounded subset $Y\subseteq\LL(\G,\mF)$, where $\G$ and $\mF$ are Banach $\wt{\C}$-modules, is bounded with respect to the ultra-pseudo-normed topology of $\LL(\G,\mF)$.

For the convenience of the reader and the sake of completeness we add a small appendix on bornological and ultrabornological $\wt{\C}$-modules. This serves as a brief survey on the subject and proves results needed in Section \ref{section_wilde}. 

We recall that throughout the paper when we talk of absorbent, balanced and convex subsets of a $\wt{\C}$-module we always intend the $\wt{\C}$-version of the corresponding classical definitions given in \cite[Definition 1.1]{Garetto:05a}. The paper \cite{Garetto:05a} has to be considered as the topological background of this work.

\section{De Wilde's theory for locally convex topological $\wt{\C}$-modules}
\label{section_wilde}
In this section we elaborate an adaptation of De Wilde's theory of webbed spaces \cite{DeWilde:67, DeWilde:71, DeWilde:78, Koethe:79, O:82} to the context of locally convex topological $\wt{\C}$-modules. Our aim is to find suitable classes of $\wt{\C}$-modules in this way to be employed in a formulation of the closed graph theorem for $\wt{\C}$-linear maps.

\subsection{Webs in locally convex topological $\wt{\C}$-modules and closed graph theorems}
\label{subsection_1}
\subsubsection*{Definitions and first properties}
\begin{definition}
\label{def_web}
Let $\G$ be a locally convex topological $\wt{\C}$-module. A family $\mathcal{W}=\{C(n_1,...,n_k)\}$ of subsets of $\G$ where $k$ and $n_1,...,n_k$ run through all of $\N\setminus 0$ is called a \emph{web} if it satisfies the relations
\[
\G=\displaystyle\bigcup_{n_1=1}^\infty C({n_1}),\qquad\qquad\qquad
C(n_1,...,n_{k-1})=\displaystyle\bigcup_{n_k=1}^\infty C(n_1,...,n_k),
\]
for $k>1$ and for all $n_1,...,n_k$.

If all the sets of a web are closed or absolutely convex we say that the web is closed respectively absolutely convex. 

A web $\mathcal{W}$ is of \emph{type $\mathcal{C}$} if the following condition is fulfilled: for all $n_k\in\N\setminus 0$ there exists $\rho_k\in\R$ such that for all $\lambda_k\ge\rho_k$ and $x_k\in C(n_1,...,n_k)$ the series $\sum_{k=1}^\infty [(\eps^{\lambda_k})_\eps]x_k$ converges in $\G$.

A web $\mathcal{W}$ is called \emph{strict} if it is absolutely convex and the following condition is fulfilled: for all $n_k\in\N\setminus 0$ there exists $\rho_k\in\R$ such that for all $\lambda_k\ge\rho_k$ and $x_k\in C(n_1,...,n_k)$ the series $\sum_{k=1}^\infty [(\eps^{\lambda_k})_\eps]x_k$ converges in $\G$ and $\sum_{k=k_0}^\infty[(\eps^{\lambda_k})_\eps]x_k:=\{u\in\G:\, \sum_{k=k_0}^n[(\eps^{\lambda_k})_\eps]x_k\to u\ \text{as}\ n\to +\infty\}$ is contained in $C(n_1,...,n_{k_0})$ for all $k_0$.  
\end{definition}
Clearly a strict web is an absolutely convex web of type $\mathcal{C}$. Conversely we have the following proposition.

\begin{proposition}
\label{prop_strict_web}
If $\mathcal{W}$ is an absolutely convex and closed web of type $\mathcal{C}$ on $\G$ then it is strict.
\end{proposition}
\begin{proof}
By definition of a web of type $\mathcal{C}$ we are able to choose the sequence $(\rho_k)_k$ such that $\rho_k\ge 0$ for all $k\ge 1$. Hence, for $x_k\in C(n_1,...,n_k)$, $\lambda_k\ge\rho_k$ and for all $k_0,N\ge 1$ we have that $\sum_{k=k_0}^{k_0+N}[(\eps^{\lambda_k})_\eps]x_k\in C(n_1,...,n_{k_0})$. Finally, since $C(n_1,...,n_{k_0})$ is closed and the sequence $(\sum_{k=k_0}^{n}[(\eps^{\lambda_k})_\eps]x_k)_n$ is convergent in $\G$ we have that  $\sum_{k=k_0}^{\infty}[(\eps^{\lambda_k})_\eps]x_k\subseteq C(n_1,...,n_{k_0})$.
\end{proof}

\begin{proposition}
\label{prop_Frechet}
Every Fr\'echet $\wt{\C}$-module has an absolutely convex and closed web of type $\mathcal{C}$.
\end{proposition}
\begin{proof}
Let $\G$ be a Fr\'echet $\wt{\C}$-module and $\{U_k\}_{k=1}^\infty$ be a base of closed and absolutely convex neighborhood of the origin with $U_{k+1}\subseteq U_k$. The subsets $C(n_1,n_2,...,n_k):=\cap_{j=1}^k[(\eps^{-n_j})_\eps]U_j$ define for $n_1,...,n_k$ running through $\N\setminus 0$ an absolutely convex and closed web of type $\mathcal{C}$. Indeed, $\G=\cup_{n_1=1}^\infty [(\eps^{-n_1})_\eps]U_1$, $C(n_1,...,n_{k-1})=\cup_{n_k=1}^\infty C(n_1,...,n_k)$ and every $C(n_1,...,n_k)$ is closed and absolutely convex by construction. Given the sequence $(n_k)_k$ let us choose $\rho_k\in\R$ such that $\rho_k-n_k\ge 0$. By the absolute convexity of the neighborhoods $U_k$ we have that for all $\lambda_k\ge\rho_k$ and $x_k\in C(n_1,...,n_k)$  
\[
\sum_{k=m}^{m+p}[(\eps^{\lambda_k})_\eps]x_k\in\sum_{k=m}^{m+p}[(\eps^{\lambda_k-n_k})_\eps]U_k\subseteq\sum_{k=m}^{m+p}[(\eps^{\rho_k-n_k})_\eps]U_m\subseteq U_m,\qquad m\ge1,\, p\ge 1.
\]
This shows that $\sum_{k=1}^\infty[(\eps^{\lambda_k})_\eps]x_k$ is a Cauchy series and then it is convergent in $\G$.
\end{proof}
We postpone a more detailed collection of examples and properties of webs to Subsection \ref{subsection_webs}. 

\subsubsection*{Closed graph theorem}
By regarding the classical De Wilde's approach as a blueprint, different versions of the closed graph theorem can be formulated in the context of \emph{webbed $\wt{\C}$-modules}. 
\begin{theorem}
\label{theorem_closed}
Let $\G$ be a Fr\'echet $\wt{\C}$-module and $\mF$ be a locally convex topological $\wt{\C}$-module which has an absolutely convex web of type $\mathcal{C}$. If the $\wt{\C}$-linear map $T:\G\to\mF$ has sequentially closed graph then it is continuous.
\end{theorem}
The proof of Theorem \ref{theorem_closed} requires the following topological lemma.

\begin{lemma}
\label{lemma_lc}
Let $\G$ and $\mF$ be locally convex topological $\wt{\C}$-modules and $T:\G\to\F$ be a $\wt{\C}$-linear map. If $\G$ is a Baire space then for all neighborhoods $V$ of $0$ in $\mF$ the set $\overline{T^{-1}(V)}$ is a neighborhood of $0$ in $\G$.
\end{lemma}
\begin{proof}
It is not restrictive to assume that $V$ is an absolutely convex neighborhood of $0$ in $\mF$. Since it is absorbent, too, we can write $\G$ as $\cup_{n=1}^\infty [(\eps^{-n})_\eps]T^{-1}(V)$ and from the assumption that $\G$ is a Baire $\wt{\C}$-module we deduce that there exists some $[(\eps^{-n})_\eps]T^{-1}(V)$ which is not rare. Hence $T^{-1}(V)$ is not rare, i.e., ${\rm{int}}(\overline{T^{-1}(V)})\neq\emptyset$. Let $x\in{\rm{int}}(\overline{T^{-1}(V)})$. Then, there exists a neighborhood $U$ of $0$ such that $x+U\subseteq \overline{T^{-1}(V)}$. Since $\overline{T^{-1}(V)}$ is absolutely convex we have that $-x\in\overline{T^{-1}(V)}$ and hence $U=-x+x+U\subseteq\overline{T^{-1}(V)}+\overline{T^{-1}(V)}\subseteq\overline{T^{-1}(V)}$. This means that $\overline{T^{-1}(V)}$ is a neighborhood of $0$ in $\G$.
\end{proof}


\begin{proof}[Proof of Theorem \ref{theorem_closed}]
Our aim is to prove that if $W$ is a closed and absolutely convex neighborhood of the origin in $\mF$ then $T^{-1}(W)$ is a neighborhood of $0$ in $\G$. Let $\mathcal{W}$ be an absolutely convex web of type $\mathcal{C}$ on $\mF$. Since $\mF=\cup_{n_1=1}^\infty C(n_1)$ we have that $\G=\cup_{n_1=1}^\infty T^{-1}(C(n_1))$ and $T^{-1}(C(n_1,...,n_{k-1}))=\cup_{n_k=1}^\infty T^{-1}(C(n_1,...,n_k))$. $\G$ is a Fr\'echet $\wt{\C}$-module and then a Baire $\wt{\C}$-module. It follows that there exists $n_1\in\N$ such that ${\rm{int}}\,\overline{T^{-1}(C(n_1))}\neq\emptyset$. From the fact that $T^{-1}(C(n_1))$ is not meager we obtain that for some $n_2\in\N$ the set $\overline{T^{-1}(C(n_1,n_2))}$ has nonempty interior. Thus, we find a sequence $\{T^{-1}(C(n_1,...,n_k))\}_k$ of sets which are not meager and absolutely convex. Writing $\mF$ as $\cup_{m=1}^\infty [(\eps^{-m})_\eps]W$ we get $T^{-1}(C(n_1,...,n_k))=\cup_{m=1}^\infty T^{-1}(C(n_1,...,n_k)\cap [(\eps^{-m})_\eps]W)$ and consequently we can extract a sequence $\{T^{-1}(C(n_1,...,n_k)\cap[(\eps^{-m_k})_\eps]W)\}_k$ of nonrare and absolutely convex subsets of $\G$. Given the sequence of real numbers $\rho_k$ associated to the choice of $(n_k)_k$ and to the web $\mathcal{W}$ we choose $\lambda_k\ge \rho_k$ such that $\lambda_k-m_k\ge 0$ for all $k$ and we define
\beq
\label{set_M_k}
M_k:=T^{-1}([(\eps^{\lambda_k})_\eps]C(n_1,...,n_k)\cap[(\eps^{\lambda_k-m_k})_\eps]W).
\eeq
By the previous considerations $M_k$ is a nonrare absolutely convex subset of $\G$ and $\overline{M_k}$ has an interior point. It follows that every $\overline{M_k}$ is a neighborhood of $0$ in $\G$. Let $(U_k)_k$ be a decreasing base of neighborhoods of the origin in $\G$ and $V_k:=\overline{M_k}\cap U_k$. If $x_0\in\overline{T^{-1}(W)}$ then there exists $x_1\in T^{-1}(W)$ such that $x_0-x_1\in V_1$ and since $V_1\subseteq\overline{M_1}$ there exists $x_2\in M_1$ such that $x_0-x_1-x_2\in V_2$. By iteration we find a sequence of points $x_{k+1}\in M_k$ for all $k\ge 1$ such that $ x_0-\sum_{k=1}^n x_k\in V_n\subseteq U_n$. In other words $\sum_{k=1}^\infty x_k$ converges to $x_0$ in $\G$. The definition \eqref{set_M_k} yields that $[(\eps^{-\lambda_k})_\eps]T(x_{k+1})\in C(n_1,...,n_k)\cap[(\eps^{-m_k})_\eps]W$ for all $k\ge 1$. Hence by definition of a web of type $\mathcal{C}$ the series 
\beq
\label{series_T}
\sum_{k=1}^\infty T(x_{k+1})=\sum_{k=1}^\infty[(\eps^{\lambda_k})_\eps][(\eps^{-\lambda_k})_\eps]T(x_{k+1})
\eeq
is convergent in $\mF$. The absolute convexity of $W$ leads to
\[
\sum_{k=1}^n T(x_k)=T(x_1)+\sum_{k=1}^{n-1}T(x_{k+1})\subseteq W+\sum_{k=1}^{n-1}[(\eps^{\lambda_k-m_k})_\eps]W\subseteq W
\]
for all $n\ge 1$ and then the series in \eqref{series_T} converges to a point of $W$. Finally, by the sequential closedness of the graph of $T$ we are allowed to conclude that the series $\sum_{k=1}^\infty T(x_k)$ converges in $W$ to $T(x_0)$ and therefore $x_0\in T^{-1}(W)$. 
We have proved that $T^{-1}(W)=\overline{T^{-1}(W)}$. By Lemma \ref{lemma_lc} this means that $T^{-1}(W)$ is a neighborhood of $0$ in $\G$.
\end{proof}


The general case of $\mF$ with a web of type $\mC$ is a little more complicated.
\begin{theorem}
\label{theorem_closed_1}
Let $\G$ be a Fr\'echet $\wt{\C}$-module and $\mF$ be a locally convex topological $\wt{\C}$-module which has a web of type $\mathcal{C}$. If the $\wt{\C}$-linear map $T:\G\to\mF$ has sequentially closed graph then it is continuous.
\end{theorem}
\begin{proof}
Let $W$ be a closed and absolutely convex neighborhood of the origin in $\mF$ and $(U_k)_k$ be a decreasing base of neighborhoods of the origin in $\G$. As in the proof of Theorem \ref{theorem_closed} we find a sequence $\{T^{-1}(C(n_1,...,n_k)\cap[(\eps^{-m_k})_\eps]W)\}_k$ of nonrare subsets of $\G$, we choose $\lambda_k\ge \rho_k$ such that $\lambda_k-m_k\ge 0$ for all $k$ and we define the sets $M_k$ as in \eqref{set_M_k}. Since every $M_k$ is not rare, for all $k$ we find $x_k\in {M_k}$ and an absolutely convex neighborhood $V_k$ of $0$ such that $V_k\subseteq U_k$ and $x_k+V_k\subseteq\overline{M_k}$. Let $x_0\in\overline{T^{-1}(W)}$. By induction arguments we find a sequence $(y_i)_{i\ge 1}$ with $y_1\in T^{-1}(W)$ and $y_k\in M_{k-1}$ for all $k\ge 2$, such that
\beq
\label{points_constr}
x_0-\sum_{i=1}^k y_i+\sum_{i=1}^{k-1}x_i\in V_k\subseteq U_k\qquad\text{and}\qquad x_0-\sum_{i=1}^k y_i+\sum_{i=1}^{k}x_i\in\overline{M_k}.
\eeq
By construction, $Tx_i\in[(\eps^{\lambda_i})_\eps]C(n_1,...,n_i)$ and since the web $\mathcal{W}$ is of type $\mC$ the series $\sum_{i=1}^\infty Tx_i$ is convergent in $\mF$. Analogously, since $y_{i+1}\in M_i$ and $Ty_{i+1}\in [(\eps^{\lambda_i})_\eps]C(n_1,...,n_i)$, the series $\sum_{i=1}^\infty Ty_i$ is convergent in $\mF$. By definition of the sets $M_i$ we easily see that $Ty_{i+1}\in [(\eps^{\lambda_i-m_i})_\eps]W$ and $Tx_i\in[(\eps^{\lambda_i-m_i})_\eps]W$. The closedness and absolute convexity of $W$ yields that $(\sum_{i=1}^k Ty_i-\sum_{i=1}^{k-1}Tx_i)_k$ is convergent to some $y_0\in W$. From \eqref{points_constr} we know that $\sum_{i=1}^k y_i-\sum_{i=1}^{k-1}x_i$ converges to $x_0$. This combined with the sequential closedness of the graph of $T$ means that $y_0=Tx_0$, i.e., $x_0\in T^{-1}(W)$. Since we have proved that $\overline{T^{-1}(W)}=T^{-1}(W)$ the continuity of the map $T$ follows from Lemma \ref{lemma_lc}.
\end{proof}


\begin{corollary}
\label{cor_ultrabor}
Let $\G$ and $(\G_\gamma)_{\gamma\in\Gamma}$ be locally convex topological $\wt{\C}$-modules such that every $\G_\gamma$ is a Fr\'echet $\wt{\C}$-module and $\G$ is endowed with the finest locally convex $\wt{\C}$-linear topology which makes the $\wt{\C}$-linear maps $\iota_\gamma:\G_\gamma\to\G$ continuous. Let $\mF$ be a locally convex topological $\wt{\C}$-module which has a web of type $\mathcal{C}$. If the $\wt{\C}$-linear map $T:\G\to\mF$ has sequentially closed graph then it is continuous.
\end{corollary}
\begin{proof}
It suffices to prove that any $\wt{\C}$-linear map $T\circ\iota_\gamma:\G_\gamma\to\mF$ is continuous. This is clear from Theorem \ref{theorem_closed_1} since the graph of $T\circ\iota_\gamma$ is sequentially closed.
\end{proof}

A further extension of Theorem \ref{theorem_closed_1} consists in assuming that $\G$ is a Baire locally convex topological $\wt{\C}$-module. In this context, in order to state the closed graph theorem the hypothesis of sequential closedness of the graph is replaced by the stronger assumption of closedness. In the proof of Theorem \ref{theorem_closed_Baire} we will use the fact that if in a topological $\wt{\C}$-module the series $\sum_{k=0}^\infty x_k$ is convergent then the corresponding sequence $x_k$ tends to $0$. 
\begin{theorem}
\label{theorem_closed_Baire}
Let $\G$ be a Baire locally convex topological $\wt{\C}$-module and $\mF$ be a locally convex topological $\wt{\C}$-module which has a web of type $\mathcal{C}$. If the $\wt{\C}$-linear map $T:\G\to\mF$ has a closed graph then it is continuous.
\end{theorem}
\begin{proof}
Let $W$ be a closed absolutely convex neighborhood of $0$ in $\mF$. Since $\G$ is a Baire locally convex topological $\wt{\C}$-module we may argue as in the proof of Theorem \ref{theorem_closed} and find a sequence of subsets 
\[
M_k:=T^{-1}([(\eps^{\lambda_k})_\eps]C(n_1,...,n_k)\cap[(\eps^{\lambda_k-m_k})_\eps]W),
\]
where $\lambda_k\ge \rho_k$ and $\lambda_k-m_k\ge 0$, such that ${\rm{int}}(\overline{M_k})\neq \emptyset$. This means that there exist $z_k\in \overline{M_k}$ and an absolutely convex neighborhood $U_k$ of $0$ in $\G$ such that $z_k+U_k\subseteq\overline{M_k}$. These relations combined with the notion of closure of a set imply that there exists a sequence of points $x_k\in M_k$ such that $x_k+U_k\subseteq\overline{M_k}$. Assume now that $x_0\in\overline{T^{-1}(W)}$. Our aim is to prove that $x_0\in T^{-1}(W)$. We begin by observing that there exists $y_1\in T^{-1}(W)$ such that $x_0-y_1\in U_1$ and $x_0-y_1+x_1\in\overline{M_1}$. We leave it to the reader to verify by induction that we can construct a sequence of points $y_k\in M_{k-1}$, such that for all $k\ge 1$, $x_0-\sum_{i=1}^k y_i+\sum_{i=1}^{k-1}x_i\in U_k$ and 
\beq
\label{2rel}
x_0-\sum_{i=1}^k y_i+\sum_{i=1}^{k}x_i\in \overline{M_k}.
\eeq
By definition of a web of type $\mathcal{C}$ both the series $\sum_{i=1}^\infty Ty_i$ and $\sum_{i=1}^\infty Tx_i$ are convergent in $\mF$. Hence the sequence $\sum_{i=1}^k Ty_i-\sum_{i=1}^{k-1} Tx_i$ is convergent. The absolute convexity of $W$ combined with the properties $Ty_1\in W$, $Ty_{i}\in[(\eps^{\lambda_{i-1}-m_{i-1}})_\eps]W$ for $i\ge 2$, $Tx_i\in[(\eps^{\lambda_i-m_i})_\eps]W$ for $i\ge 1$ entails $\sum_{i=1}^{k}Ty_i-\sum_{i=1}^{k-1} Tx_i\in W$ for all $k$. Since $W$ is closed this means that the sequence $\sum_{i=1}^k Ty_i-\sum_{i=1}^{k-1} Tx_i$ converges to some $y_0\in W$. We complete the proof by showing that $(x_0,y_0)\in\overline{{\rm{Graph}}(T)}$. Indeed, if this holds the assumption of closedness of the graph of $T$ implies $y_0=T(x_0)$ and $x_0\in T^{-1}(W)$. Let $U,V$ be neighborhoods of $0$ in $\G$ and $\mF$ respectively. By \eqref{2rel} it follows that there exists $t_k\in M_k$ such that $x_0-\sum_{i=1}^k y_i+\sum_{i=1}^{k}x_i-t_k\in U$. Since the series $\sum_{k=1}^\infty Tt_k$ and $\sum_{k=1}^\infty Tx_k$ are convergent we have that $Tt_k\to 0$ and $Tx_k\to 0$. This combined with $\sum_{i=1}^k Ty_i-\sum_{i=1}^{k-1} Tx_i\to y_0$ leads to $y_0-\sum_{i=1}^k Ty_i+\sum_{i=1}^{k}Tx_i-Tt_k\to 0$. Hence, $y_0-\sum_{i=1}^k Ty_i+\sum_{i=1}^{k}Tx_i-Tt_k\in V$ for all $k$ larger than some $k_0$. In conclusion, 
\[
(x_0,y_0)-\biggl(\sum_{i=1}^k y_i-\sum_{i=1}^{k}x_i+t_k,T\biggl(\sum_{i=1}^k y_i-\sum_{i=1}^{k}x_i+t_k\biggr) \biggr)\in (U,V)
\]
for $k\ge k_0$, so that $(x_0,y_0)\in\overline{{\rm{Graph}}(T)}$.
\end{proof}

We focus now our attention on \emph{ultrabornological} $\wt{\C}$-modules in order to provide some further extensions of Theorem \ref{theorem_closed_1} and Corollary \ref{cor_ultrabor}.  

We begin by observing that when $\G$ is a locally convex topological $\wt{\C}$-module and $A$ is a \emph{disk}, i.e., a nonempty, balanced and convex subset of $\G$, then the $\wt{\C}$-submodule generated by $A$ is simply $\G_A:=\cup_{n\in\N}[(\eps^{-n})_\eps]A$. On $\G_A$ the function $\val_A(u):=\sup\{b\in\R:\, u\in [(\eps^b)_\eps]A\}$ is a valuation and $\mP_A(u):=\esp^{-\val_A(u)}$ is an ultra-pseudo-seminorm. If $A$ is absorbent then $\G_A$ coincides with the whole of $\G$.

In coherence with the usual language adopted for locally convex topological vector spaces, we say that $A\subseteq\G$ is a \emph{Banach disk} if it is a disk of $\G$ and $(\G_A,\mP_A)$ is a Banach $\wt{\C}$-module. 
\begin{definition}
\label{def_born}
A subset $A$ of a locally convex topological $\wt{\C}$-module $\G$ is said to be \emph{ultrabornivorous} if it absorbs any bounded Banach disk of $\G$.\\ 
A locally convex topological $\wt{\C}$-module $\G$ is \emph{ultrabornological} if every absolutely convex and ultrabornivorous subset of $\G$ is a neighborhood of the origin.
\end{definition}
Every Banach $\wt{\C}$-module $(\G,\mQ)$ is ultrabornological. Indeed, taking the absorbent, absolutely convex and bounded subset $A:=\{u\in\G:\, \mQ(u)\le 1\}$, we have that $\G=\G_A$ and $\mP_A=\mQ$. It follows that if $V$ is an absolutely convex and ultrabornivorous subset of $\G$ then it absorbs $\G_A$, that is $\G_A\subseteq [(\eps^b)_\eps]V$ for all $b\le a$. Hence, $[(\eps^{-a})_\eps]A\subseteq V$ implies that $V$ is a neighborhood of $0$ in $\G$.

It is possible to characterize the topology of a Hausdorff ultrabornological $\wt{\C}$-module as the finest locally convex $\wt{\C}$-linear topology which makes any injection $\iota_B:\G_B\to \G$ continuous, where $B$ is any bounded Banach disk. This is proved in Proposition \ref{prop_born_ultra} of the appendix at the end of the paper. The appendix provides a topological investigation of the class of ultrabornological $\wt{\C}$-modules and of the wider class of bornological $\wt{\C}$-modules. In particular it proves that every Fr\'echet $\wt{\C}$-module is ultrabornological.

It is now clear that Corollary \ref{cor_ultrabor} applies to a Hausdorff ultrabornological $\wt{\C}$-module $\G$.
\begin{proposition}
\label{prop_cor_ultrabor}
Let $\G$ be a Hausdorff ultrabornological $\wt{\C}$-module and $\mF$ be a locally convex topological $\wt{\C}$-module which has a web of type $\mathcal{C}$. If the $\wt{\C}$-linear map $T:\G\to\mF$ has sequentially closed graph then it is continuous.
\end{proposition}

\subsection{Webbed $\wt{\C}$-modules: examples and hereditary properties}
\label{subsection_webs}

\begin{proposition}
\label{prop_seq_closed}
Every sequentially closed $\wt{\C}$-submodule $\mH$ of a locally convex topological $\wt{\C}$-module $\G$ equipped with a web $\mathcal{W}$ of type $\mathcal{C}$, has a web $\mathcal{V}$ of type $\mathcal{C}$. Moreover, if $\mathcal{W}$ is absolutely convex then $\mathcal{V}$ is absolutely convex.
\end{proposition}
\begin{proof}
Let $\mathcal{W}=\{C(n_1,...,n_k)\}$. It is immediate to prove that a web of type $\mC$ on $\mH$ is given by $\mathcal{V}:=\{C(n_1,...,n_k)\cap \mH\}$. Indeed, if $(\rho_k)_k$ is the sequence associated to $(C(n_1,...,n_k))_k$, $\lambda_k\ge \rho_k$ and $x_k\in C(n_1,...,n_k)\cap\mH$ then by the sequentially closedness of $\mH$ we have that the series $\sum_{k=1}^\infty [(\eps^{\lambda_k})_\eps]x_k$ is convergent in $\mH$.
\end{proof}
The following result is straightforward from the notion of sequential continuity.
\begin{proposition}
\label{prop_image}
Let $\G$ and $\mF$ be locally convex topological $\wt{\C}$-modules and $T:\G\to\mF$ be a $\wt{\C}$-linear sequentially continuous map. If $\mW=\{C(n_1,...,n_k)\}$ is a web of type $\mathcal{C}$ on $\G$ then $T(\mW)=\{T(C(n_1,...,n_k))\}$ is a web of type $\mC$ on $T(\G)$. If $\mW$ is absolutely convex then $T(\mW)$ is absolutely convex.
\end{proposition}
\begin{proposition}
\label{prop_quotient}
Let $\G$ be a locally convex topological $\wt{\C}$-module, $M$ a $\wt{\C}$-submodule of $\G$ and $\mW$ a web of type $\mC$ on $\G$. The quotient $\G/M$ equipped with the quotient topology has a web $\mV$ of type $\mC$. Moreover, if $\mW$ is absolutely convex then $\mV$ is absolutely convex.
\end{proposition}
\begin{proof}
We recall that as observed in \cite[Example 1.12]{Garetto:05a} $\G/M$ equipped with the quotient topology is a locally convex topological $\wt{\C}$-module and the $\wt{\C}$-linear map $\pi:\G\to\G/M:u\to u+M$ is continuous. Hence, by Proposition \ref{prop_image} the assertion is clear.
\end{proof}
We investigate now the properties, as a webbed $\wt{\C}$-module, of the dual $\LL(\G,\wt{\C})$ of a strict inductive limit of locally convex topological $\wt{\C}$-modules. We begin with the following result on strict inductive limits.
\begin{proposition}
\label{prop_strict}
Let $\G$ be the strict inductive limit of the sequence of locally convex topological $\wt{\C}$-modules $(\G_p)_p$. If on any $\G_p$ there exists a web $\mathcal{W}_p$ of type $\mathcal{C}$ then $\G$ has a web $\mW$ of type $\mC$. Moreover, if every $\mW_p$ is absolutely convex then $\mW$ is absolutely convex.
\end{proposition}
\begin{proof}
Given $\mW_p=\{C_p(n_1,...,n_k)\}$ we define $D(n_1):=\G_{n_1}$ and $D(n_1,...,n_k):=C_{n_1}(n_2,...,n_k)$ for $k\ge 2$. The collection of sets $\{D(n_1,...,n_k)\}$ is a web on $\G$. Indeed, by construction it is clear that $\cup_{n_1=1}^\infty D(n_1)=\G$ and $D(n_1,...,n_{k-1})=\cup_{n_k=1}^\infty D(n_1,...,n_{k-1},n_k)$. Let $(n_k)_{k\ge 1}$ be a sequence in $\N\setminus 0$. Consider the web $\mW_{n_1}$ on $\G_{n_1}$ and take the sequence $(\rho_k)_{k\ge 2}$ associated to $(n_k)_{k\ge 2}$ (the choice of $\rho_1$ is arbitrary). When $\lambda_k\ge\rho_k$ and $x_k\in D(n_1,...,n_k)$ we have that $x_k\in C_{n_1}(n_2,...,n_k)\subseteq\G_{n_1}$. Since $\mW_{n_1}$ is a web of type $\mC$ on $\G_{n_1}$ we have that the series $\sum_{k=2}^\infty[(\eps^{\lambda_k})_\eps]x_k$ is convergent in $\G_{n_1}$. Hence the series $\sum_{k=1}^\infty[(\eps^{\lambda_k})_\eps]x_k$ is convergent in $\G$.
\end{proof}
We proceed with a lemma and a proposition useful for the proof of Theorem \ref{theorem_dual} concerning the dual $\LL(\G,\wt{\C})$.
\begin{lemma}
\label{lemma_C}
A web $\mW=\{C(n_1,...,n_k)\}$ on a locally convex topological $\wt{\C}$-module $\G$ is of type $\mC$ if the following condition is satisfied: for all sequences $(n_k)_{k\ge 1}$ there exists $(\mu_k)_{k\ge 1}\subseteq\R$ such that the sequence $([(\eps^{\mu_k})_\eps]x_k)_{k\ge 1}$ with $x_k\in C(n_1,...,n_k)$ is contained in a bounded, absolutely convex and sequentially complete subset $M$ of $\G$.
\end{lemma}
\begin{proof}
Let $(n_k)_k$ a sequence of natural numbers and $x_k\in C(n_1,...,n_k)$. Choosing $\rho_k$ such that $\rho_k-\mu_k\ge k$ we have that for $\lambda_k\ge\rho_k$,
\[
\sum_{k=1}^N[(\eps^{\lambda_k})_\eps]x_k=\sum_{k=1}^N[(\eps^{\lambda_k-\mu_k})_\eps][(\eps^{\mu_k})_\eps]x_k\in\sum_{k=1}^N[(\eps^{\lambda_k-\mu_k})_\eps]M\subseteq M. 
\]
The sequence $(\sum_{k=1}^N[(\eps^{\lambda_k})_\eps]x_k)_k$ is a Cauchy sequence in $M$. Indeed, the boundedness of $M$ implies that for any balanced neighborhood $U$ of $0$ in $\G$ there exists $a\in\R$ such that $M\subseteq[(\eps^a)_\eps] U$ and when $N+1+a\ge 0$ we get that
\[
\sum_{k=1}^{N+p}[(\eps^{\lambda_k})_\eps]x_k-\sum_{k=1}^{N}[(\eps^{\lambda_k})_\eps]x_k\in [(\eps^{N+1})_\eps]M\subseteq[(\eps^{N+1+a})_\eps] U\subseteq U.
\]
Since $M$ is sequentially complete this means that the series $\sum_{k=1}^\infty [(\eps^{\lambda_k})_\eps]x_k$ is convergent in $M$.
\end{proof}
\begin{proposition}
\label{prop_top}
Let $\G$ be a topological $\wt{\C}$-module.
\begin{trivlist}
\item[(i)] If $A$ is a neighborhood of $0$ in $\G$ then $A^\circ:=\{T\in\LL(\G,\wt{\C}):\, \forall u\in A\ |Tu|_\esp\le 1\}$ is equicontinuous.
\item[(ii)] If $D\subseteq\LL(\G,\wt{\C})$ is equicontinuous then it is $\beta_b(\LL(\G,\wt{\C}),\G)$-bounded.
\end{trivlist}
\end{proposition}
\begin{proof}
$(i)$ Combining the fact that $A$ is a neighborhood of $0$ in $\G$ with the definition of $A^\circ$ we have that for any $\eta>0$ on the neighborhood $U:=[(\eps^a)_\eps]A$ with $\esp^{-a}\le\eta$ the inequality $|T(u)|_\esp\le\eta$ holds for all $u\in U$ and $T\in A^\circ$. This proves that $A^\circ$ is equicontinuous.

$(ii)$ Let $B$ be a bounded subset of $\G$. Then for all neighborhoods $U$ of $0$ in $\G$ there exists $b\in\R$ such that $[(\eps^{-b})_\eps]B\subseteq U$. Since $D$ is equicontinuous there exists  a neighborhood $U$ of $0$ such that $|T(u)|_\esp\le 1$ for all $u\in U$ and $T\in D$. This means that for all $T\in D$,
\[
\sup_{u\in B}|T(u)|_\esp=\esp^{-b}|T([(\eps^{-b})_\eps]u)|_\esp\le \esp^{-b}.
\]
The previous estimate proves that $D$ is a $\beta_b(\LL(\G,\wt{\C}),\G)$-bounded subset of the dual $\LL(\G,\wt{\C})$.
\end{proof}
\begin{theorem}
\label{theorem_dual}
Let $\G$ be the strict inductive limit of the sequence of metrizable $\wt{\C}$-modules $(\G_p)_p$. Then the dual $\LL(\G,\wt{\C})$ endowed with the topology $\beta_b(\LL(\G,\wt{\C}),\G)$ has an absolutely convex and closed web of type $\mathcal{C}$.
\end{theorem}
\begin{proof}
Let $U_1^p\supseteq U_2^p\supseteq U_3^p\supseteq...$ be a countable base of neighborhoods of the origin in $\G_p$. We set
\[
C(n_1,...,n_k):=\bigcap_{j=1}^k(U^j_{n_j})^{\circ}.
\]
We begin by proving that for all $p\in\N$, $\LL(\G,\wt{\C})=\cup_{m=1}^\infty (U^p_m)^{\circ}$.  Indeed, if $T\in\LL(\G,\wt{\C})$ then $T|_{\G_p}\in\LL(\G_p,\wt{\C})$ and for some neighborhood $U^p_m$ we have that $|T(u)|_\esp\le 1$ for all $u\in U^p_m$. As a consequence $\LL(\G,\wt{\C})=\cup_{n_1=1}^\infty C(n_1)$ and $C(n_1,...,n_{k-1})=\cup_{n_k=1}^\infty C(n_1,...,n_k)$. By \cite[Proposition 2.4]{Garetto:05a}
it follows that $\mW=\{C(n_1,...,n_k)\}$ is an absolutely convex and closed web on $\LL(\G,\wt{\C})$. It remains to prove that it is of type $\mC$. Let $(n_k)_{k\ge 1}$ be a sequence in $\N\setminus 0$ and $T_{k}\in C(n_1,...,n_k)$. If we show that $(T_{k})_k$ is contained in a bounded, absolutely convex and sequentially complete subset $M$ of $\LL(\G,\wt{\C})$ endowed with the topology $\beta_b(\LL(\G,\wt{\C}),\G)$ then Lemma \ref{lemma_C} yields that $\mW$ is of type $\mC$. First we note that for all $p\in\N$ there exists $m_p\ge 1$ such that $(T_k)_{k\ge 1}\subseteq (U_{m_p}^p)^\circ$. In fact, for $k\ge p$ we have by construction that $T_k\in C(n_1,...,n_k)=\cap_{j=1}^k(U_{n_j}^j)^\circ\subseteq\cap_{j=1}^p(U_{n_j}^j)^\circ\subseteq (U^p_{n_p})^\circ$ and since $\LL(\G,\wt{\C})=\cup_{l=1}^\infty (U^p_l)^\circ$ there exists $l_p$ such that $T_k\in (U^p_{l_p})^\circ$ for all $k\le p$. Hence, taking $m_p=\max\{n_p,l_p\}$  we get $(T_k)_{k\ge 1}\subseteq\cap_{p\in\N}(U^p_{m_p})^\circ$. The set $M:=\cap_{p\in\N}(U^p_{m_p})^\circ$ is absolutely convex and $\sigma(\LL(\G,\wt{\C}),\G)$-closed, so $\beta_b(\LL(\G,\wt{\C}),\G)$-closed. Since the dual $\LL(\G,\wt{\C})$ endowed with the topology of uniform convergence on bounded subsets is complete, we conclude that $M$ is complete itself. Let $A$ be the absolutely convex hull of $\cup_{p\in\N}U^p_{m_p}$. Since it is the set of all finite $\wt{\C}$-linear combinations $\sum_{i=1}^n\lambda_i x_i$ with $x_i\in U^i_{m_i}$ and $\max_{i=1,...,n}|\lambda_i|_\esp\le 1$ one can easily prove $A$ is a neighborhood of $0$ in $\G$ and that $M= A^\circ$. By Proposition \ref{prop_top} it follows that $M$ is equicontinuous and therefore $\beta_b(\LL(\G,\wt{\C}),\G)$-bounded. Concluding, $M$ is an absolutely convex, complete and bounded subset of $(\LL(\G,\wt{\C}),\beta_b(\LL(\G,\wt{\C}),\G))$ and contains the sequence $(T_k)_k$.
\end{proof}
Further hereditary properties, here omitted for the sake of brevity, concern topological products and projective limits of locally convex topological $\wt{\C}$-modules equipped with webs of type $\mC$, and are easily obtained by adapting the corresponding classical arguments in \cite[Ch. 4, Sec. 4]{DeWilde:78} and \cite[Ch. 35, Sec. 4]{Koethe:79}.

\section{Strict morphisms of $\wt{\C}$-modules and open mapping theorems for locally convex or metrizable topological $\wt{\C}$-modules}
\label{section_open}
A first kind of open mapping theorems can be obtained via application of the closed graph theorems stated in the previous section. This will require the choice of suitable classes of locally convex topological $\wt{\C}$-modules as webbed and ultrabornological $\wt{\C}$-modules.


\begin{theorem}
\label{theorem_open}
Let $\G$ be a locally convex topological $\wt{\C}$-module $\G$ with a web of type $\mC$ and $\mF$ be a Hausdorff ultrabornological $\wt{\C}$-module. If $T:\G\to\mF$ is a $\wt{\C}$-linear, surjective and continuous map then it is open.
\end{theorem}
\begin{proof}
By Proposition \ref{prop_quotient} we know that $\G/{\rm{ker}}\, T$ has a web of type $\mC$. Since $T$ is continuous and $\G/{\rm{ker}}\, T$ is endowed with the quotient topology, the $\wt{\C}$-linear map $\overline{T}:\G/{\rm{ker}}\, T\to\mF:u+{\rm{ker}}\, T\to T(u)$ is continuous and then $\overline{T}^{-1}$ has a sequentially closed graph. By applying Proposition \ref{prop_cor_ultrabor} to $\overline{T}^{-1}$ we conclude that $\overline{T}^{-1}$ is continuous. Then for all open subsets $A$ of $\G$ the image $\overline{T}(A+{\rm{ker}}\, T)=T(A)$ is an open subset of $\mF$, which means that the map $T$ is open.
\end{proof} 
An interesting extension of the previous theorem consists in replacing the assumption of continuity of $T$ with the assumption of having a sequentially closed graph. A map which has sequentially closed graph is said to be sequentially closed.


\begin{theorem}
\label{theorem_open_closed_1}
Let $\G$ be a locally convex topological $\wt{\C}$-module with a web of type $\mC$ and $\mF$ be a Fr\'echet $\wt{\C}$-module. Let $T$ be a $\wt{\C}$-linear map defined on a $\wt{\C}$-submodule $D$ of $\G$ such that $T(D)$ is not meager. If ${\rm{Graph}}(T)$ is sequentially closed in $\G\times\mF$ then $T(D)=\mF$ and $T$ is open.
\end{theorem}
\begin{proof}
Let $\mathcal{W}=\{C(n_1,...,n_k)\}$ be a web of type $\mC$ on $\G$. Since $T(D)=\bigcup_{n_1=1}^{+\infty}T(C(n_1)\cap D)$ we find $n_1$ such that $\text{int}(\overline{T(C(n_1)\cap D)})\neq\emptyset$. Recursively we can define a sequence $n_1,n_2,...$ such that the set $T(C(n_1,...,n_k)\cap D)$ is not meager in $\mF$. Let $V$ be an absolutely convex and closed neighborhood of the origin in $\G$. Then for all $k$ there exists $m_k\in\N$ such that $T(C(n_1,...,n_k)\cap D\cap[(\eps^{-m_k})_\eps]V)$ is not rare. By definition of a web of type $\mC$ we have a sequence of real numbers $\rho_k$ such that for all $\lambda_k\ge\rho_k$ and for all $x_k\in C(n_1,...,n_k)$ the series $\sum_k [(\eps^{\lambda_k})_\eps]x_k$ is convergent. Let us fix $\lambda_k\ge\rho_k$ such that $\lambda_k-m_k\ge 0$ for all $k$ and define
\[
M_k:=[(\eps^{\lambda_k})_\eps](C(n_1,...,n_k)\cap D\cap[(\eps^{-m_k})_\eps]V).
\]
By construction $T(M_k)$ is not rare in $\mF$. Let $U_k$ be a decreasing sequence of neighborhoods of the origin in $\mF$. Since $T(M_k)$ is not rare we find $x_k\in M_k$ and an absolutely convex neighborhood $U^{(k)}\subseteq U_k$ of the origin such that $Tx_k+U^{(k)}\subseteq \overline{T(M_k)}$. Moreover, since $T(D)=\cup_{n=1}^\infty [(\eps^{-n})_\eps]T(D\cap V)$ is not meager then $\overline{T(D\cap V)}$ contains a neighborhood of the origin in $\mF$.

Let now $y_0$ be a point of $\overline{T(D\cap V)}$. Then there exists $z_1\in D\cap V$ such that $y_0-Tz_1\in U^{(1)}$. Consequently, $y_0-Tz_1+Tx_1\in\overline{T(M_1)}$. By induction we find, for every $i\ge 2$, an element $z_i$ in $M_{i-1}$ such that 
\[
y_0-\sum_{1}^k Tz_i+\sum_{1}^{k-1}Tx_i\in U^{(k)}\subseteq U_k
\]
and
\[
y_0-\sum_{1}^k Tz_i+\sum_{1}^{k}Tx_i\in \overline{T(M_k)},
\]
for all $k\ge 1$. Clearly $\big(\sum_{1}^k Tz_i-\sum_{1}^{k-1}Tx_i\big)_k$ converges to $y_0$ in $\mF$. Since $z_i\in[(\eps^{\lambda_{i-1}})_\eps]C(n_1,...,n_{i-1})$ for $i\ge 2$, the sequence $(\sum_{1}^k z_i)_k$ is convergent. Analogously the sequence $(\sum_i^{k-1}x_i)_k$ is convergent in $\G$. It follows that $\sum_{1}^k z_i-\sum_i^{k-1}x_i\in D$ converges in $\G$ to some point $x_0\in V$. Using the fact that the graph of $T$ is sequentially closed in $\G\times\mF$ we conclude that $x_0\in D$ and $Tx_0=y_0$. This means that $y_0\in T(D\cap V)$ and therefore $T(D\cap V)=\overline{T(D\cap V)}$ is a neighborhood of the origin in $\mF$. It is immediate to deduce that $T(D)=\mF$ and that the map $T$ is open.
\end{proof}


\begin{theorem}
\label{theorem_open_closed_2}
Let $\G$ be a locally convex topological $\wt{\C}$-module with a web of type $\mC$ and $\mF$ be a Hausdorff ultrabornological $\wt{\C}$-module. A sequentially closed $\wt{\C}$-linear map $T$ of $\G$ onto $\mF$ is open.
\end{theorem}
\begin{proof}
Let $B$ be a bounded Banach disk of $\mF$. We denote the restriction of $T$ to $T^{-1}(\mF_B)$ by $T_B:T^{-1}(\mF_B)\to\mF_B$. $\mF_B$ is a Banach $\wt{\C}$-module, $T^{-1}(\mF_B)$ is a $\wt{\C}$-submodule of $\G$ whose image is not meager and $T_B$ has a sequentially closed graph. Hence, by Theorem \ref{theorem_open_closed_1} we conclude that $T_B$ is an open map. Let $V$ be an absolutely convex neighborhood of the origin in $\G$. Since $V_B:=V\cap T^{-1}(\mF_B)$ is a neighborhood of the origin in $T^{-1}(\mF_B)$ we obtain that $U_B\subseteq T_B(V_B)$ where $U_B$ is a neighborhood of $0$ in $\mF_B$. Finally, the fact that $V$ is absolutely convex and $T_B(V_B)\subseteq T(V)$ implies that $T(V)$ contains the absolutely convex hull $W$ of $\cup_B U_B$. But by Proposition \ref{prop_born_ultra}$(iii)$ $W$ is a neighborhood of $0$ in $\mF$ and therefore $T(V)$ is a neighborhood of $0$ in $\mF$. 
\end{proof}


We now consider the wider family of topological $\wt{\C}$-modules. Our aim is providing (cf. Theorem \ref{theorem_hor_1}) an open mapping theorem valid for topological $\wt{\C}$-modules which are not necessarily locally convex but metrizable. Inspired by \cite{Horvath:66} we introduce the following definitions.
\begin{definition}
\label{def_strict_mor_1}
A bijective continuous $\wt{\C}$-linear map $T$ from a topological $\wt{\C}$-module $\G$ onto a topological $\wt{\C}$-module $\mF$ is called an \emph{isomorphism} if the inverse map $T^{-1}$ is continuous. An injective continuous $\wt{\C}$-linear map $T:\G\to\mF$ is a \emph{strict morphism} if it is an isomorphism from $\G$ onto $T(\G)$. 
\end{definition}
\begin{definition}
\label{def_strict_mor_2}
A continuous $\wt{\C}$-linear map $T$ from a topological $\wt{\C}$-module $\G$ into a topological $\wt{\C}$-module $\mF$ is a \emph{strict morphism} if the associated injection $\overline{T}:\G/{\rm{Ker}}(T)\to\mF$ is a strict morphism. 
\end{definition}
\begin{proposition}
\label{prop_eq_strict_mor}
If $T$ is a continuous $\wt{\C}$-linear map from a topological $\wt{\C}$-module $\G$ into a topological $\wt{\C}$-module $\mF$ then the following conditions are equivalent:
\begin{itemize}
\item[(i)] $T$ is a strict morphism;
\item[(ii)] $T$ maps every neighborhood of $0$ in $\G$ onto a neighborhood of $0$ in $T(\G)$;
\item[(iii)] $T$ maps every open set of $\G$ onto an open set of $T(\G)$.
\end{itemize}
\end{proposition}
We omit the proof of Proposition \ref{prop_eq_strict_mor} since it follows the lines of argument of the analogous result valid in the topological vector space context.
%
%
%
The proof of the open mapping theorem \ref{theorem_hor_1} requires two preliminary lemmas.  
\begin{lemma}
\label{lemma_hor_1}
Let $\G$ and $\mF$ be two topological $\wt{\C}$-modules and let $T:\G\to\mF$ be a continuous surjective $\wt{\C}$-linear map. If $\mF$ is a Baire $\wt{\C}$-module, then for every neighborhood $U$ of $0$ in $\G$ the set $\overline{T(U)}$ is a neighborhood of $0$ in $\mF$.
\end{lemma} 
\begin{proof}
Proposition 1.3 in \cite{Garetto:05a} allows to find a balanced and absorbent neighborhood $V$ of 0 such that $V+V\subseteq U$. It follows that we can write $\G$ as $\cup_{n\in\N}[(\eps^{-n})_\eps]V$ and $\mF=T(\G)=\cup_{n\in\N}[(\eps^{-n})_\eps]T(V)$. $\mF$ is a Baire space. Therefore, there exists $n\in\N$ such that $\overline{[(\eps^{-n})_\eps]T(V)}=[(\eps^{-n})_\eps]\overline{T(V)}$ has non empty interior. Since $T(V)$ as well as its closure are balanced subsets of $\mF$ we conclude that $0$ belongs to the interior of $\overline{T(V)}+\overline{T(V)}$. The continuity of the addition in $\mF$ yields $\overline{T(V)}+\overline{T(V)}\subseteq\overline{T(V)+T(V)}$. But $T(V)+T(V)=T(V+V)\subseteq T(U)$ and hence $\overline{T(V)}+\overline{T(V)}\subseteq \overline{T(U)}$. This means that $0$ belongs to the interior of $\overline{T(U)}$ or in other words that $\overline{T(U)}$ is a neighborhood of $0$ in $\mF$.
\end{proof}

The following lemma is reported in \cite[Chapter 3, Section 17, Lemma 2]{Horvath:66}.
\begin{lemma}
\label{lemma_hor_2}
Let $\G$ and $\mF$ be two metric spaces and assume that $\G$ is complete. Suppose that $T$ is a continuous map from $\G$ into $\mF$ which has the following property: for every $r>0$ there exists $\rho>0$ such that for all $u\in\G$ the image $T(B_r(u))$ of the set
\[
B_r(u):=\{v	\in\G:\quad d_\G(u,v)\le r\}
\]
is dense in the ball
\[
B_\rho(T(u)):=\{w\in\mF:\quad d_{\mF}(w,T(u))\le\rho\}.
\]
Then for every $a>r$ the set $T(B_a(u))$ contains $B_\rho(T(u))$.
\end{lemma}


\begin{theorem}
\label{theorem_hor_1}
Let $\G$ and $\mF$ be two metrizable complete topological $\wt{\C}$-modules and $T$ a continuous surjective $\wt{\C}$-linear map from $\G$ onto $\mF$. Then $T$ is a strict morphism.
\end{theorem}
\begin{proof}
$\G$ and $\mF$ are both Baire $\wt{\C}$-modules. Then Lemma \ref{lemma_hor_1} implies that $\overline{T(B_r(0))}$ is a neighborhood of $0$ in $\mF$. Hence there exists $\rho>0$ such that $B_\rho(0)\subseteq\overline{T(B_r(0))}$, i.e., $T(B_r(0))$ is dense in $B_\rho(0)\subseteq\mF$. Since we can assume that the distances on $\G$ and $\mF$ are invariant under translation (the proof of this fact is totally analogous to \cite[Ch. 2, Sec 6, Th. 1]{Horvath:66}) we obtain that $B_\rho(T(u))\subseteq\overline{T(B_r(0))}+T(u)=\overline{T(B_r(u))}$. This means that for every $u\in\G$ the set $T(B_r(u))$ is dense in $B_\rho(T(u))$. By Lemma \ref{lemma_hor_2} we conclude that $B_\rho(T(u))\subseteq T(B_a(u))$ for all $a>r$. In particular $B_\rho(0)\subseteq T(B_a(0))$ and therefore $T(B_a(0))$ is a neighborhood of $0$ in $\mF$. Since $r$ is arbitrary we have that the second assertion of Proposition \ref{prop_eq_strict_mor} is satisfied or equivalently that $T$ is a strict morphism.
\end{proof}


Theorem \ref{theorem_hor_1} yields the following version of the closed graph theorem.
\begin{theorem}
\label{theorem_hor_2}
Let $\G$ and $\mF$ be two metrizable complete topological $\wt{\C}$-modules and $T:\G\to\mF$ a $\wt{\C}$-linear map whose graph is closed in $\G\times\mF$. Then $T$ is continuous.
\end{theorem}
\begin{proof}
We begin by recalling that the product of two metrizable complete topological $\wt{\C}$-modules is complete when endowed with the product topology. By assumption Graph$(T)$ is a closed $\wt{\C}$-submodule of $\G\times\mF$ and so it is metrizable and complete. The projection $\pi_1:{\rm{Graph}}(T)\to\G:(u,T(u))\to u$ is bijective, $\wt{\C}$-linear and continuous and therefore by Theorem \ref{theorem_hor_1} its inverse is continuous. If $\pi_2$ is the projection of Graph$(T)$ into $\mF$ we conclude that $T=\pi_2\circ\pi_1^{-1}$ is continuous.
\end{proof}
\begin{corollary}
\label{corollary_hor_2}
Let $\G$ and $\mF$ be two metrizable complete topological $\wt{\C}$-modules and $T:\G\to\mF$ a $\wt{\C}$-linear map fulfilling the following property: for every sequence of points $(u_n)_n$ tending to $0$ in $\G$ and for which $(Tu_n)_n$ tends to some $v\in\mF$, one has necessarily $v=0$. Then $T$ is continuous.
\end{corollary}
\begin{proof}
Suppose that $(u,w)$ adheres to the graph of $T$. Then there exists a sequence of points $(u_n)_n$ such that $u_n\to u$ and $Tu_n\to w$. But then $u_n-u\to 0$ and $T(u_n-u)\to w-T(u)$. By our assumption it follows that $w=T(u)$. This means that Graph$(T)$ is closed and thus by Theorem \ref{theorem_hor_2} the map $T$ is continuous. 
\end{proof}

\section{Strings and $\wt{\C}$-linear topologies without convexity conditions}
\label{section_adasch}
At the end of Section \ref{section_open} we have proved a version of the closed graph theorem which holds for topological $\wt{\C}$-modules which are not necessarily locally convex but metrizable. In this section it is our intention to provide closed graph and open mapping theorems for a larger class of topological $\wt{\C}$-modules. For such a purpose we elaborate a suitable adaptation of the theory of strings for topological vector spaces (see \cite{Adasch:78}) to the setting of topological $\wt{\C}$-modules. 

\subsection{Strings: definition and basic properties}
\label{subsec_strings}
\begin{definition}
\label{def_string}
Let $\G$ be a $\wt{\C}$-module. A sequence $\mathcal{U}=(U_n)_{n\in\N\setminus 0}$ of subsets of $\G$ is a \emph{string} in $\G$ if
\begin{trivlist}
\item[(i)] every $U_n$ is balanced,
\item[(ii)] every $U_n$ is absorbent,
\item[(iii)] $U_{n+1}+U_{n+1}\subseteq U_n$ for all $n$,
\item[(iv)] $[(\eps^{-n})_\eps]U_{n+p}\subseteq U_p$ for all $n,p$.
\end{trivlist}
\end{definition}
$U_1$ is called the beginning of the string $\mathcal{U}$ while $U_n$ is the n-th knot of $\mathcal{U}$. If $\mU=(U_n)_n$ and $\mV=(V_n)_n$ are strings in $\G$ and $\lambda\in\wt{\C}$ we define $\lambda \mU:=(\lambda U_n)_n$, the sum $\mU+\mV:=(U_n+V_n)_n$ and the intersection $\mU\cap\mV:=(U_n\cap V_n)_n$. When $\lambda$ is invertible these are all strings of subsets of $\G$. 

Given a string $\mU$ the set $N(\mU):=\cap_n U_n$ is called the \emph{kernel} of $\mU$. $N(\mU)$ is a $\wt{\C}$-submodule of $\G$. The additivity follows from $(iii)$. Assume now that $u\in N(\mU)$ and take $\lambda\in\wt{\C}$ with $\val(\lambda)<0$ (the fact that $u\in N(\mU)$ implies $\lambda u\in N(\mU)$ is clear when $\val(\lambda)\ge 0$). Let us choose $n\in\N$ such that $\val(\lambda)+n\ge 0$. For this $n$ let us fix $p\in\N$ such that $p-n\ge 1$. Combining $(i)$ with $(iv)$ we obtain that if $u\in N(\mU)$ then 
\[
\lambda u=[(\eps^{+\val({\lambda})})_\eps][(\eps^{-\val({\lambda})})_\eps]\lambda u\in [(\eps^{+\val({\lambda})})_\eps]U_p\subseteq[(\eps^{+\val({\lambda})+n})_\eps][(\eps^{-n})_\eps]U_p\subseteq [(\eps^{-n})_\eps]U_p\subseteq U_{p-n}
\]  
Varying $p$ we conclude that $\lambda u\in N(\mU)$.

Finally for $\mU=(U_n)_n$ and $\mV=(V_n)_n$ we use the notation $\mU\subseteq \mV$ when $U_n\subseteq V_n$ for all indices $n$.

\begin{definition}
\label{def_top_string}
Let $\G$ be a topological $\wt{\C}$-module. A string $\mU$ in $\G$ is said to be topological if every knot is a neighborhood of the origin in $\G$.
\end{definition}
\begin{remark}
\label{remark_top_modules}
Every neighborhood of the origin $U$ in a topological $\wt{\C}$-modules $\G$ generates a topological string. We begin by taking a balanced and absorbent neighborhood $U_1$ of $0$ which is contained in $U$. Since the assumption $(iv)$ in Definition \ref{def_string} is equivalent to require that $U_{n+1}\subseteq [(\eps)_\eps]U_n$ for all $n$, by iteration we can choose a sequence of balanced neighborhoods $U_n$ of the origin such that $U_{n+1}\subseteq [(\eps)_\eps]U_n$ and $U_{n+1}+U_{n+1}\subseteq U_n$. This gives a topological string $(U_n)_n$ in $\G$.
\end{remark}

\begin{proposition}
\label{prop_string_1}
Let $\G$ be a topological $\wt{\C}$-module. Then there exists a set $\mathcal{S}$ of strings in $\G$ such that:
\begin{trivlist}
\item[(i)] if $\mU\in\mathcal{S}$ and $\mV\in\mS$ then there exists $\mW\in\mS$ such that $\mW\subseteq \mU\cap\mV$,
\item[(ii)] the knots of the strings in $\mS$ form a base of neighborhoods of $0$ in $\G$.
\end{trivlist}
Moreover, the topology on $\G$ is Hausdorff if and only if $\displaystyle\cap_{\mU\in\mS} N(\mU)=\{0\}$.
\end{proposition}
\begin{proof}
Let $\mathcal{B}$ be a base of neighborhoods of the origin in $\G$. By Remark \ref{remark_top_modules} for each $B\in\mathcal{B}$ we can construct a topological string $(U_{B,n})_n$ with $U_{B,1}\subseteq B$. Let $\mS$ be the set of all finite intersections of strings constructed in this way. Clearly, properties $(i)$ and $(ii)$ are fulfilled. The final assertion is a simple reformulation of the separatedness of $\G$. 
\end{proof}
It is clear that the set of all topological strings of a topological $\wt{\C}$-module $\G$ satisfies the conditions $(i)$ and $(ii)$ of Proposition \ref{prop_string_1}.

\begin{proposition}
\label{prop_string_2}
Let $\mS$ be a set of strings in a $\wt{\C}$-module $\G$ such that the following condition is fulfilled:

for all $\mU\in\mS$ and $\mV\in\mS$ there exists $\mW\in\mS$ such that $\mW\subseteq \mU\cap\mV$.

Then the knots of the strings in $\mS$ form a base of neighborhoods of the origin for a $\wt{\C}$-linear topology $\tau_\mS$ on $\G$. Furthermore, if $\cap_{\mU\in\mS} N(\mU)=\{0\}$ then $\tau_\mS$ is Hausdorff.
\end{proposition}
\begin{proof}
Let $\mB$ be the set of all the knots of the strings in $\mS$. We define a topology $\tau_\mS$ on $\G$ by claiming that $U(x)$ is a neighborhood of $x\in\G$ if and only if there exists $U\in\mB$ such that $x+U\subseteq U(x)$. We want to prove that $\tau_\mS$ is $\wt{\C}$-linear and that the knots of the strings of $\mS$ form a base of neighborhoods of the origin. The addition $+:\G\times\G\to \G$ is continuous with respect to the topology $\tau_\mS$. Indeed, let $W({u_0+v_0})$ be a neighborhood of $u_0+v_0\in\G$. Then, there exists a knot $W_{n_0}$ of a certain string $(W_n)_n$ of $\mS$ such that $u_0+v_0+W_{n_0}\subseteq W(u_0+v_0)$. By the third defining property of a string we have that $u_0+W_{n_0+1}+v_0+W_{n_0+1}\subseteq W(u_0+v_0)$ where $u_0+W_{n_0+1}$ and $v_0+W_{n_0+1}$ are neighborhoods of $u_0$ and $v_0$ respectively. We now want to prove that the scalar product between generalized numbers and elements of $\G$ is continuous. Let $\lambda_0\in\wt{\C}$, $u_0\in\G$ and $U(\lambda_0u_0)$ be a neighborhood of $\lambda_0u_0$. By definition of $\tau_\mS$ there exists a knot $U_{n_0}$ of a certain string $(U_n)_n\in\mS$ such that $\lambda_0 u_0+U_{n_0}\subseteq U(\lambda_0u_0)$. Let us write $\lambda u-\lambda_0 u_0$ as $\lambda(u-u_0)+(\lambda-\lambda_0)u_0$. Since the knots of a string are absorbent subsets there exists $b\in\R$ such that $u_0\in [(\eps^b)_\eps]U_{n_0+1}$ and assuming that $\val(\lambda-\lambda_0)\ge -b$ we have that $(\lambda-\lambda_0)u_0\in U_{n_0+1}$. It follows that $\val(\lambda)\ge\min\{-b,\val(\lambda_0)\}$. Let $n\in\N$ such that $\min\{-b,\val(\lambda_0)\}\ge -n$. Taking $u-u_0\in U_{n_0+1+n}$ by property $(iv)$ in Definition \ref{def_string} we conclude that if $\val(\lambda-\lambda_0)\ge -b$ then $\lambda(u-u_0)\in U_{n_0+1}$. Therefore, $\lambda u-\lambda_0 u_0\in U_{n_0+1}+U_{n_0+1}\subseteq U_{n_0}$ when $\val(\lambda-\lambda_0)\ge -b$ and $u\in u_0+U_{n_0+1+n}$ with $\min\{-b,\val(\lambda_0)\}\ge -n$.
\end{proof}
A set $\mS$ of strings in $\G$ with the property  $(i)$ of Proposition \ref{prop_string_1} is called \emph{directed}. A directed set of strings $\mS$ in a topological $\wt{\C}$-module $(\G,\tau)$ such that $\tau_\mS=\tau$ is called \emph{fundamental}.

Since the set of all strings in a $\wt{\C}$-module $\G$ is directed from the previous proposition we have that it generates a $\wt{\C}$-linear topology $\tau^f$ on $\G$. By Proposition \ref{prop_string_1} this is the finest $\wt{\C}$-linear topology on $\G$. 

Every absolutely convex and absorbent subset $U$ of a $\wt{\C}$-module $\G$ gives a string $\mU_U=(U_n)_n$ with $U_n=[(\eps^{n})_\eps]U$. Indeed, $U_{n+1}+U_{n+1}=[(\eps^{n+1})_\eps](U+U)\subseteq(\eps^{n})_\eps]U=U_n$ and $U_{n+1}=[(\eps)_\eps]U_n$. $\mU_U$ is called the \emph{natural string of $U$}.
\begin{proposition}
\label{prop_string_convex}
Let $\mS_c$ be the set of the natural strings of all the absorbent and absolutely convex subsets of $\G$. Then,  $\tau^c:=\tau_{\mS_c}$ is a locally convex $\wt{\C}$-linear topology on $\G$. 
\end{proposition}
\begin{proof}
Since the set $\mS_c$ is directed, by Proposition \ref{prop_string_2} we have that $\tau^c$ is a $\wt{\C}$-linear topology with a base of absorbent and absolutely convex neighborhoods of the origin. Hence it is a locally convex $\wt{\C}$-linear topology.
\end{proof}
$\tau^c$ is the finest locally convex $\wt{\C}$-linear topology on $\G$ and by construction $\tau^c$ is coarser than $\tau^f$.


\subsection{Barrelled topological $\wt{\C}$-modules}
\label{subsec_barrelled}
In the sequel we denote the categories of topological $\wt{\C}$-modules and locally convex topological $\wt{\C}$-modules by $\mathcal{L}$ and $\mC$ respectively. In addition, we say that a string in a topological $\wt{\C}$-module $\G$ is closed if every knot is closed.
\begin{definition}
\label{def_barrelled}
A topological $\wt{\C}$-module $\G$ is $\mL$-barrelled if every closed string in $\G$ is topological.
\end{definition}
Note that a notion of barrelledness (or more precisely $\mC$-barrelledness) already exists within the category of locally convex topological $\wt{\C}$-modules and it is given by Definition 2.12 in \cite{Garetto:05a}. A locally convex topological $\wt{\C}$-module $\G$ which is $\mL$-barrelled is also $\mC$-barrelled. Indeed, if $U$ is a barrel of $\G$, i.e. it is an absolutely convex, absorbent and closed subset of $\G$, then the string $\mU_U$ is natural and closed and consequently by Definition \ref{def_barrelled} every knot of $\mU_U$ is a neighborhood of the origin. This means that $U$ is a neighborhood of the origin in $\G$. Since we only deal with $\mL$-barrelledness without comparison with $\mC$-barrelledness we allow ourselves to omit the prefix $\mL$ in this paper. 

\begin{proposition}
\label{prop_Baire}
Every Baire topological $\wt{\C}$-module is barrelled.
\end{proposition}
\begin{proof}
Let $\mU=(U_n)_n$ be a closed string. For fixed $n$ we have that $\G=\cup_{k\in\N}[(\eps^{-k})_\eps]U_{n+1}$. Since $\G$ is a Baire space it follows that $U_{n+1}$ has an interior point, i.e. there exists $u\in\G$ and an open balanced neighborhood $V$ of $0$ such that $u+V\subseteq U_{n+1}$. By the third property characterizing the knots of a string we obtain that $u+V-u+V\subseteq U_{n+1}+U_{n+1}\subseteq U_n$ and therefore $U_n$ is a neighborhood of $0$ in $\G$.
\end{proof}
Since every Fr\'echet $\wt{\C}$-module is a Baire $\wt{\C}$-module we have from the previous proposition that every Fr\'echet $\wt{\C}$-module is barrelled. 

\begin{proposition}
\label{prop_barrelled_quo}
The quotient of a barrelled $\wt{\C}$-module $\G$ with respect to a $\wt{\C}$-submodule $M$ is a barrelled $\wt{\C}$-module when endowed with the quotient topology.
\end{proposition}
\begin{proof}
We denote the original topology on $\G$ and the quotient topology on $\G/M$ by $\tau$ and $\widehat{\tau}$ respectively. $\pi:\G\to\G/M$ is the canonical projection of $\G$ onto $\G/M$. It is clear that $(\G/M,\widehat{\tau})$ is a topological $\wt{\C}$-module. Let $\mU=(U_n)_n$ be a closed string in $(\G/M,\widehat{\tau})$. Then, $\pi^{-1}(\mU):=(\pi^{-1}(U_n))_n$ is a string in $\G$ and by definition of the quotient topology each $\pi^{-1}(U_n)$ is closed. Since $\G$ is barrelled it follows that $\pi^{-1}(U_n)$ is a neighborhood of $0$ in $\G$ and therefore $\mU$ is a topological string in the quotient $\G/M$.
\end{proof}

\begin{definition}
\label{def_tau_b}
Let $\mS_b$ the set of all closed strings in a topological $\wt{\C}$-module $\G$. We call \emph{strong topology} on $\G$ the topology $\tau^b=\tau_{\mS_b}$ generated by $\mS_b$.
\end{definition}
If $\tau$ is the original $\wt{\C}$-linear topology on $\G$ then $\tau\prec\tau^b$, where $\prec$ stands for ``coarser than''. This is due to the fact that in a topological $\wt{\C}$-module every neighborhood of $0$ contains a closed and balanced neighborhood of $0$.

\begin{proposition}
\label{prop_barrelled_1}
Let $\G$ be a topological $\wt{\C}$-module endowed with the topology $\tau$. The following conditions are equivalent:
\begin{trivlist}
\item[(i)] $\G$ is barrelled.
\item[(ii)] $\tau=\tau^b$.
\item[(iii)] Every $\wt{\C}$-linear topology $\tau'$ on $\G$ which has a base of $\tau$-closed neighborhoods of the origin is coarser than $\tau$.
\end{trivlist}
\end{proposition}
\begin{proof}
$(i)\Rightarrow (ii)$ We already know that $\tau\prec\tau^b$. Let $U$ be a neighborhood of $0$ for $\tau^b$. Then there exists a closed string $(U_n)_n$ such that $U_n\subseteq U$ for some $n$. Since $\G$ is barrelled we have that $U_n$ is a neighborhood of $0$ for the topology $\tau$ and hence the same property holds for $U$.

$(ii)\Rightarrow (i)$ Let $\mU=(U_n)_n$ be a closed string in $(\G,\tau)$. Every $U_n$ is a neighborhood of $0$ with respect to $\tau^b$ and by $(ii)$ it follows that every $U_n$ is a neighborhood of $0$ for the topology $\tau$. This means that $\G$ is barrelled.

$(ii)\Rightarrow (iii)$ Let $U$ be a neighborhood of the origin for the topology $\tau'$. By definition there exists a $\tau$-closed neighborhood $V$ of $0$ such that $V\subseteq U$. Since $V$ generates a closed topological string it follows that $U$ is a neighborhood of $0$ with respect to $\tau^b$. By $(ii)$ we can now conclude that $U$ is a neighborhood of $0$ for $\tau$ and therefore $\tau'\prec\tau$.

$(iii)\Rightarrow (ii)$ Since $\tau^b$ is a $\wt{\C}$-linear topology which has a base of $\tau$-closed neighborhoods of the origin the assertion $(iii)$ yields $\tau^b\prec\tau$. Thus, $\tau^b=\tau$.
\end{proof}


\begin{proposition}
\label{prop_barrelled_cont}
Let $(\G,\tau_1)$ and $(\mF,\tau_2)$ be topological $\wt{\C}$-modules. If the $\wt{\C}$-linear map $T:(\G,\tau_1)\to(\mF,\tau_2)$ is continuous then $T:(\G,\tau_1^b)\to(\mF,\tau_2^b)$ is continuous.
\end{proposition}
\begin{proof}
The continuity and the $\wt{\C}$-linearity of $T$ yield that if $\mU=(U_n)_n$ is a closed string in $\mF$ then $T^{-1}(\mU)=(T^{-1}U_n)_n$ is a closed string in $\G$. Let now $U$ be a neighborhood of $0$ in $(\mF,\tau_2^b)$. By definition of $\tau_2^b$ there exists a closed string $\mU=(U_n)_n$ and some $n\in\N$ such that $U_n\subseteq U$. Hence, $T^{-1}(U_n)\subseteq T^{-1}(U)$. From the considerations above we have that $T^{-1}(U)$ is a neighborhood of $0$ with respect to $\tau_1^b$.
\end{proof}


We now investigate the barrelledness properties of inductive limits of topological $\wt{\C}$-modules. First of all we recall that given a family $(\G_i,\tau_i)_{i\in I}$ of topological $\wt{\C}$-modules and a family of $\wt{\C}$-linear maps $f_i:\G_i\to \G$, where $\G$ is the $\wt{\C}$-span of the images $f_i(\G_i)$, we say that $\G$ is the inductive limit of the $\wt{\C}$-modules $(\G_i)_{i\in I}$ if it is endowed with the finest $\wt{\C}$-linear topology which makes every map $f_i$ continuous. One can easily prove that the set of strings
\[
\mS_{\text{i.l.}}:=\{\mU:\ \mU\, \text{is a string in $\G$ and $f_i^{-1}(\mU)$ is a topological string in $\G_i$ for all $i$}\}
\]
is directed and that the inductive limit topology on $\G$ coincides with $\tau_{\mS_{\rm{i.l.}}}$.\\
When the $\wt{\C}$-modules $\G_i$ are barrelled the following hereditary property holds.
\begin{proposition}
\label{prop_induc_barrelled}
The inductive limit $\G$ of a family of barrelled $\wt{\C}$-modules $(\G_i,\tau_i,f_i)_{i\in I}$ is barrelled.
\end{proposition}
\begin{proof}
Let $\mU$ be a closed string in $\G$. Then $f_i^{-1}(\mU)$ is a closed string in $\G_i$ and since every $\G_i$ is barrelled $f_i^{-1}(\mU)\subseteq \G_i$ is a topological string for all $i\in I$. This means that $\mU\in\mS_{\rm{i.l.}}$ and therefore $\mU$ is a topological string for the inductive limit topology on $\G$.
\end{proof}

A topological $\wt{\C}$-module $(\G,\tau)$ endowed with the finest topology $\tau^f$ is barrelled. Indeed, $\tau^f\prec(\tau^f)^b$ and by definition of the finest topology $(\tau^b)^f\prec\tau^f$. Proposition \ref{prop_barrelled_1} yields that $(\G,\tau^f)$ is barrelled. We can now look for the coarsest barrelled $\wt{\C}$-linear topology which is finer then $\tau$. This is obtained as the inductive limit topology determined by the family of barrelled $\wt{\C}$-linear topologies finer then $\tau$. By the previous proposition this topology is itself barrelled.
\begin{definition}
\label{def_tau_t}
We call the coarsest barrelled $\wt{\C}$-linear topology on $(\G,\tau)$ which is finer than $\tau$ the \emph{associated barrelled topology}. It is denoted by $\tau^t$.
\end{definition}

\begin{proposition}
\label{prop_barrelled_cont_t}
Let $(\G,\tau_1)$ and $(\mF,\tau_2)$ be topological $\wt{\C}$-modules. If the $\wt{\C}$-linear map $T:(\G,\tau_1)\to(\mF,\tau_2)$ is continuous then $T:(\G,\tau_1^t)\to(\mF,\tau_2^t)$ is continuous.
\end{proposition}
\begin{proof}
We take the finest $\wt{\C}$-linear topology on $\mF$ such that the map $T:(\G,\tau_1^t)\to\mF$ is continuous and denote this topology by $\tau'$. Proposition \ref{prop_barrelled_cont} entails the continuity of the map $T:(\G,(\tau_1^t)^b)\to(\mF,(\tau')^b)$, where $(\tau_1^t)^b=\tau_1^t$. As a consequence $(\tau')^b\prec\tau'$ and by Proposition \ref{prop_barrelled_1} we conclude that $(\mF,\tau')$ is barrelled. It is clear that the map $T:(\G,\tau_1^t)\to(\mF,\tau_2)$ is continuous and therefore $\tau_2\prec \tau'$. This yields $\tau_2^t\prec\tau'$ and the continuity of $T:(\G,\tau_1^t)\to(\mF,\tau_2^t)$ as a straightforward consequence.
\end{proof}

\subsection{Barrelled topological $\wt{\C}$-modules and the closed graph theorem}
\label{subsec_closed_bar}
This subsection is devoted to prove a version of the closed graph theorem for barrelled $\wt{\C}$-modules. We recall that the adjective closed, when referring to a map, means that the corresponding graph is closed.
\begin{theorem}
\label{theo_closed_barrelled}
Every closed $\wt{\C}$-linear map from a barrelled $\wt{\C}$-module $(\G,\tau_1)$ to a Fr\'echet $\wt{\C}$-module $(\mF,\tau_2)$ is continuous.
\end{theorem}
The proof of Theorem \ref{theo_closed_barrelled} requires some preliminary results.

 
\begin{proposition}
\label{prop_pre_1}
Let $(\G,\tau_1)$ be a topological $\wt{\C}$-module and $(\mF,\tau_2)$ a Hausdorff topological $\wt{\C}$-module. If $T:(\G,\tau_1)\to(\mF,\tau_2)$ is a $\wt{\C}$-linear continuous map then its graph is closed with respect to the product topology on $\G\times\mF$.
\end{proposition}
\begin{proof}
Assume that $(u,v)\not\in{\rm{Graph}}(T)$, which means that $v\neq T(u)$. Then there exist neighborhoods $V$ and $U$ of $v$ and $T(u)$ respectively such that $V\cap U=\emptyset$. Since the map $T$ is continuous $T^{-1}(U)\times V$ is a neighborhood of $(u,v)$ in the product topology on $\G\times\mF$ and by construction $(T^{-1}(U)\times V)\cap {\rm{Graph}}(T)=\emptyset$.
\end{proof}


\begin{lemma}
\label{lemma_pre_2}
Let $(\G,\tau_1)$ and $(\mF,\tau_2)$ be topological $\wt{\C}$-modules where $\tau_2$ is a Hausdorff $\wt{\C}$-linear topology. The following assertions are equivalent:
\begin{itemize}
\item[(i)] the $\wt{\C}$-linear map $T:(\G,\tau_1)\to(\mF,\tau_2)$ is closed;
\item[(ii)] there exists a Hausdorff $\wt{\C}$-linear topology $\tau_3$ on $\mF$ with $\tau_3\prec\tau_2$ such that $T:(\G,\tau_1)\to(\mF,\tau_3)$ is continuous;
\item[(iii)] the strings $T(\mU)+\mV$, where $\mU$ is any topological string in $(\G,\tau_1)$ and $\mV$ is any topological string in $(\mF,\tau_2)$, generate a Hausdorff $\wt{\C}$-linear topology on $\mF$.
\end{itemize}
\end{lemma} 
\begin{proof}
$(ii)\Rightarrow(i)$ Since $T:(\G,\tau_1)\to(\mF,\tau_3)$ is continuous we have by Proposition \ref{prop_pre_1} that Graph$(T)$ is closed in $(\G,\tau_1)\times(\mF,\tau_3)$. The product topology of $(\G,\tau_1)\times(\mF,\tau_2)$ is finer than $\tau_1\times\tau_3$ and therefore the graph of $T$ is closed in $(\G,\tau_1)\times(\mF,\tau_2)$.

$(i)\Rightarrow(iii)$ Assume that Graph$(T)$ is closed with respect to the product topology on $(\G,\tau_1)\times(\mF,\tau_2)$. Then $\cap_{U,V}T(U)+V=\{0\}$ where $U$ and $V$ vary in the families of neighborhoods of $0$ in $\G$ and $\mF$, respectively. Indeed, if $y\neq 0$ then $(0,y)\not\in{\rm{Graph}}(T)$ and there exist some balanced neighborhoods $U$ and $V$ of $0$ in $\G$ and $\mF$ respectively such that $((0,y)+(U,V))\cap{\rm{Graph}}(T)=\emptyset$. It follows that $y\not\in T(U)+V$ and that the topology generated by the strings $T(\mU)+\mV$ is Hausdorff.

$(iii)\Rightarrow(ii)$ We denote the topology generated by the the strings $T(\mU)+\mV$ by $\tau_3$. By hypothesis it is Hausdorff. Moreover $\tau_3$ is coarser than $\tau_2$ and makes the map $T:(\G,\tau_1)\to(\mF,\tau_3)$ continuous. 
\end{proof}


\begin{proposition}
\label{prop_pre_Fre}
On a Fr\'echet $\wt{\C}$-module there exists no coarser Hausdorff barrelled $\wt{\C}$-linear topology.
\end{proposition}
\begin{proof}
We denote the original topology on a Fr\'echet $\wt{\C}$-module $\G$ by $\tau_0$ and we assume that there exists $\tau\prec\tau_0$ Hausdorff and barrelled. Let $(U_n)_n$ be a closed string in $(\G,\tau_0)$ whose knots form a base of neighborhoods of the origin for the topology $\tau_0$. We want to prove that this string is $\tau$-topological. This will imply that $\tau_0\prec\tau$. 

Let us consider $u\in\overline{U_{n+1}}^{\,\tau}$. Hence $(u+\overline{U_{n+2}}^{\,\tau})\cap U_{n+1}\neq\emptyset$, i.e. there exist $u_1\in U_{n+1}$ such that $u-u_1\in\overline{U_{n+2}}^{\,\tau}$. This is due to the fact that $(\overline{U_n}^{\,\tau})_n$ is a topological string since $\tau$ is a barrelled $\wt{\C}$-linear topology. In the same way we find $u_2\in U_{n+2}$ such that $u-u_1-u_2\in\overline{U_{n+3}}^{\,\tau}$ and by iteration a sequence of elements $u_j\in U_{n+j}$ such that $u-\sum_{k=1}^j u_k\in\overline{U_{n+1+j}}^{\,\tau}$. The Cauchy sequence $\big(\sum_{k=1}^j u_k\big)$ converges to some $y$ with respect to the topology $\tau_0$ and since 
\[
\sum_{k=1}^j u_k\in \sum_{k=1}^j U_{n+k}\subseteq U_n, 
\]
where $U_n$ is $\tau_0$-closed, we conclude that $y\in U_n$. Clearly $\big(\sum_{k=1}^j u_k\big)\to y\in U_n$ in the coarser topology $\tau$. 

Finally from $u-\sum_{k=1}^j u_k\in\overline{U_{n+1+j}}^{\,\tau}$ for all $j$ we obtain that $u-y\in\bigcap_{j=1}^\infty \overline{U_j}^{\,\tau}$. But $\bigcap_{j=1}^\infty \overline{U_j}^{\,\tau}=\{0\}$. Indeed, if $z\neq 0$ then there exists a neighborhood $V$ of $0$ for the topology $\tau$ such that $z+V\cap V=\emptyset$. Hence $(z+V)\cap U_j=\emptyset$ for some $j$ and therefore $z\not\in\overline{U_j}^{\, \tau}$. We have proved that $u=y\in U_n$ or in other words that $\overline{U_{n+1}}^{\,\tau}\subseteq U_n$. It follows that $U_n$ is a neighborhood of $0$
for the topology $\tau$. 
\end{proof}

\begin{proof}[Proof of Theorem \ref{theo_closed_barrelled}]
Since the map $T:(\G,\tau_1)\to(\mF,\tau_2)$ is closed, by Lemma \ref{lemma_pre_2}$(ii)$ there exists a Hausdorff topology $\tau_3\prec\tau_2$ such that the $\wt{\C}$-linear map $T:(\G,\tau_1)\to(\mF,\tau_3)$ is continuous. Proposition \ref{prop_barrelled_cont_t} yields the continuity of $T:(\G,\tau_1^t)\to(\mF,\tau_3^t)$. Since $\tau_3^t\prec\tau_2^t$ and $\tau_2^t=\tau_2$ we have that $\tau_3^t$ is a Hausdorff barrelled $\wt{\C}$-linear topology coarser than $\tau_2$. By Proposition \ref{prop_pre_Fre} we can conclude that $\tau_3^t=\tau_2$. Finally, from the fact that $(\G,\tau_1)$ is barrelled and therefore $\tau_1^t=\tau_1$ we obtain that the map $T$ is continuous from $(\G,\tau_1)$ to $(\mF,\tau_2)$.
\end{proof}

We now try to find a larger class of topological $\wt{\C}$-modules $(\mF,\tau_2)$ for which Theorem \ref{theo_closed_barrelled} can be stated.  

Let $T$ be a closed $\wt{\C}$-linear map from a barrelled topological $\wt{\C}$-module $(\G,\tau_1)$ to a Hausdorff topological $\wt{\C}$-module $(\mF,\tau_2)$. By Lemma \ref{lemma_pre_2}$(ii)$ there exists a Hausdorff $\wt{\C}$-linear topology $\tau_3\prec \tau_2$ on $\mF$ such that the map $T:(\G,\tau_1)\to(\mF,\tau_3)$ is continuous. Then by Proposition \ref{prop_barrelled_cont_t} the map $T$ is continuous from $(\G,\tau_1^t)$ to $(\mF,\tau_3^t)$. Since $\tau_1^t=\tau_1$ we are able to conclude that the map $T:(\G,\tau_1)\to(\mF,\tau_2)$ is continuous if we know that $\tau_2\prec\tau_3^t$. This motivates the following definition.
\begin{definition}
\label{def_infra_s}
A Hausdorff topological $\wt{\C}$-module $(\mF,\tau_0)$ is called an \emph{infra-s} $\wt{\C}$-module if for every coarser Hausdorff $\wt{\C}$-linear topology $\tau$ on $\mF$ we have that $\tau_0\prec \tau^t$ or equivalently $\tau^t=\tau_0^t$.  
\end{definition} 
It is immediate to extend Theorem \ref{theo_closed_barrelled} as follows.
\begin{theorem}
\label{theo_closed_barrelled_2}
Every closed $\wt{\C}$-linear map from a barrelled $\wt{\C}$-module $(\G,\tau_1)$ to an infra-s $\wt{\C}$-module $(\mF,\tau_2)$ is continuous.
\end{theorem}
\begin{example}
\label{ex_infra_s}
\leavevmode
\begin{trivlist}
\item[(i)] A first example of infra-s $\wt{\C}$-modules is given by Fr\'echet $\wt{\C}$-modules. Indeed, if $(\mF,\tau_0)$ is a Fr\'echet $\wt{\C}$-module and $\tau$ is a coarser Hausdorff $\wt{\C}$-linear topology then $\tau^t$ is a Hausdorff barrelled $\wt{\C}$-linear topology on $\mF$ which is coarser than $\tau_0^t=\tau_0$. Hence by Proposition \ref{prop_pre_Fre} $\tau^t=\tau_0$.
\item[(ii)] Let $(\mF,\tau_0)$ be a Fr\'echet $\wt{\C}$-module such that the weak topology $\sigma(\mF,\mL(\mF,\wt{\C}))$ is Hausdorff. Since $\sigma\prec\tau_0$, $\sigma^t$ is a Hausdorff barrelled topology coarser than $\tau_0$ and then by Proposition \ref{prop_pre_Fre} it coincides with $\tau_0$. Let now $\tau$ be a Hausdorff $\wt{\C}$-linear topology coarser than $\sigma$. It follows that $\tau^t\prec\sigma^t=\tau_0$ and therefore $\tau^t=\tau_0$ by Proposition \ref{prop_pre_Fre} again. Finally, $\sigma\prec\tau_0=\tau^t$ yields that $(\mF,\sigma)$ is an infra-s $\wt{\C}$-module.
\end{trivlist}
\end{example}

The statement of Theorem \ref{theo_closed_barrelled_2} provides a characterization of infra-s $\wt{\C}$-modules as we will see in the sequel.
\begin{proposition}
\label{prop_charac_infra_s}
Let $(\mF,\tau_2)$ be a Hausdorff topological $\wt{\C}$-module with the property that every closed $\wt{\C}$-linear map from a barrelled $\wt{\C}$-module $(\G,\tau_1)$ to $(\mF,\tau_2)$ is continuous. Then $(\mF,\tau_2)$ in an infra-s $\wt{\C}$-module.
\end{proposition}
\begin{proof}
Let $\tau$ be a Hausdorff $\wt{\C}$-linear topology on $\mF$ coarser than $\tau_2$. The identity map $I$ from $(\mF,\tau)$ onto $(\mF,\tau_2)$ is closed. Since $\tau\prec\tau^t$, $I$ remains closed as a map from $(\mF,\tau^t)$ onto $(\mF,\tau_2)$. $(\mF,\tau^t)$ is a barrelled $\wt{\C}$-module, then by hypothesis $I:(\mF,\tau^t)\to(\mF,\tau_2)$ is continuous. It follows that $\tau_2$ is coarser than $\tau^t$ and therefore $(\mF,\tau_2)$ is an infra-s $\wt{\C}$-module. 
\end{proof}

\begin{remark}
\label{rem_charac_infra}
\leavevmode
\begin{trivlist}
\item[(i)] The proof of the previous proposition shows that a Hausdorff topological $\wt{\C}$-module $(\mF,\tau_2)$ is an infra-s $\wt{\C}$-module if and only if every closed $\wt{\C}$-linear bijection from a barrelled $\wt{\C}$-module $(\G,\tau_1)$ onto $(\mF,\tau_2)$ is continuous.
\item[(ii)] The previous characterization as well as the statement of Proposition \ref{prop_charac_infra_s} hold under the assumption that $(\G,\tau_1)$ is Hausdorff and barrelled.
\end{trivlist}

\end{remark}

As already seen in Section \ref{section_open} in the case of locally convex topological $\wt{\C}$-modules appropriate versions of the open mapping theorem can be obtained via application of the closed graph theorem.

\subsection{Barrelled topological $\wt{\C}$-modules and the open mapping theorem}
\label{subsec_open_bar}
For the sequel of the paper it will be useful to recall the following result concerning quotients of topological $\wt{\C}$-modules.
\begin{proposition}
\label{prop_quotient_com}
Let $\G$ be a complete metrizable topological $\wt{\C}$-modules and $M$ a closed $\wt{\C}$-submodule of $\G$. Then $\G/M$ equipped with the quotient topology is complete.
\end{proposition}
We omit the proof of Proposition \ref{prop_quotient_com} since it makes use of the arguments employed in proving the analogous statement for topological vector spaces in \cite[Th. 2, Ch. 2, Sec. 9]{Horvath:66}. Combining this result with Remark 1.4 and Example 1.12 in \cite{Garetto:05a} we have that the quotient of a Fr\'echet $\wt{\C}$-module with respect to a closed $\wt{\C}$-submodule is a Fr\'echet $\wt{\C}$-module when equipped with the quotient topology. 


\begin{theorem}
\label{theorem_open_mapping_bar}
Let $\G$ be a Fr\'echet $\wt{\C}$-module and $\mF$ a Hausdorff barrelled topological $\wt{\C}$-module. If $T:\G\to\mF$ is a $\wt{\C}$-linear, surjective and continuous map then it is open.
\end{theorem}
\begin{proof}
As observed above the quotient $\G/{\rm{ker}}\, T$ has the topological structure of a Fr\'echet $\wt{\C}$-module.  Since $T$ is continuous by definition of the quotient topology the $\wt{\C}$-linear map $\overline{T}:\G/{\rm{ker}}\, T\to\mF:u+{\rm{ker}}\, T\to T(u)$ is continuous and then $\overline{T}^{-1}$ has a closed graph. This is obtained by applying Proposition \ref{prop_pre_1} to $\overline{T}$ and to the fact that Graph$(\overline{T}^{-1})=\{(v,u):\, (u,v)\in{\rm{Graph}}(\overline{T})\}$. $\overline{T}^{-1}$ is a closed $\wt{\C}$-linear map from a Hausdorff barrelled $\wt{\C}$-module to a Fr\'echet $\wt{\C}$-module and therefore by Theorem \ref{theo_closed_barrelled} it is continuous. This allows to conclude that the map $T$ is open.  
\end{proof}

In order to extend the validity of Theorem \ref{theorem_open_mapping_bar} to a larger class of topological $\wt{\C}$-modules $(\G,\tau_1)$ we continue the investigation of the properties of strict morphisms started in Section \ref{section_open}.
\begin{proposition}
\label{prop_strict_mor_ker}
Let $T$ be a continuous $\wt{\C}$-linear map from a topological $\wt{\C}$-module $\G$ into a topological $\wt{\C}$-module $\mF$. Let $\mV_0$ be the set of all the neighborhoods of $0$ in $\mF$. If $T$ is a strict morphism then  
\[
\bigcap_{V\in\mV_0}\overline{T^{-1}(V)}=\overline{\rm{Ker}(T)}.
\]
\end{proposition}
\begin{proof}
Let $\mU_0$ be the set of all neighborhoods of $0$ in $\G$. By Proposition \ref{prop_eq_strict_mor} we know that if $T$ is a strict morphism and $U\in\mU_0$ then $T(U)$ is a neighborhood of the origin in $T(\G)$ with respect to the topology induced by $\mF$. It is obvious that $\overline{\rm{Ker}(T)}\subseteq\cap_{V\in\mV_0}\overline{T^{-1}(V)}$ since $\cap_{V\in\mV_0}\overline{T^{-1}(V)}$ is closed and ${\rm{Ker}(T)}\subseteq{T^{-1}(V)}$ for all $V$. Conversely,
\[
\bigcap_{V\in\mV_0}\overline{T^{-1}(V)}=\bigcap_{U\in\mU_0,V\in\mV_0}{T^{-1}(V)+U}\subseteq \bigcap_{U\in\mU_0}{T^{-1}(T(U))+U}=\bigcap_{U\in\mU_0}(U+{\rm{Ker}(T)}+U)=\overline{\rm{Ker}(T)}.
\]
\end{proof}

\begin{definition}
\label{def_weak_sing}
Let $T$ be a $\wt{\C}$-linear map between the topological $\wt{\C}$-modules $(\G,\tau_1)$ and $(\mF,\tau_2)$. We set 
\[
{\rm{K}}(T):=\bigcap_{V\in\mV_0}\overline{T^{-1}(V)}\qquad\qquad\quad \text{and}\qquad\qquad\quad {\rm{S}}(T):=\bigcap_{U\in\mU_0}\overline{T(U)}.
\]
The set ${\rm{S}}(T)$ is called \emph{singularity} of $T$. If ${\rm{K}}(T)=\overline{\rm{Ker}(T)}$ then $T$ is called \emph{weakly singular}.  
\end{definition}
By Proposition \ref{prop_strict_mor_ker} we know that every strict morphism between topological $\wt{\C}$-modules is weakly singular. 
\begin{proposition}
\label{prop_S(T)}
${\rm{S}}(T)$ is a closed $\wt{\C}$-submodule of $\mF$. Moreover, ${\rm{S}}(T)=\{v\in\mF:\ (0,v)\in\overline{{\rm{Graph}}(T)}\}$.
\end{proposition}
\begin{proof}
Let $U\in\mU_0$ and $W\in\mU_0$. We take $U'$ and $W'$ in $\mU_0$ such that $U'+U'\subseteq U$ and $W'+W'\subseteq W$. If $u_1,u_2\in{\rm{S}}(T)$ then $(u_1+W')\cap T(U')\neq \emptyset$ and $(u_2+W')\cap T(U')\neq \emptyset$. It follows that $(u_1+u_2+W)\cap T(U)\neq\emptyset$. Since the choice of the neighborhoods $U$ and $W$ was arbitrary, we obtain that $u_1+u_2\in{\rm{S}}(T)$. 

Let $u\in{\rm{S}}(T)$ and $\lambda\in\wt{\C}$. By the continuity of the scalar multiplication map $\G\to\G:v\to\lambda v$ we know that for all $U,W\in\mU_0$ we can find $U',W'\in\mU_0$ such that $\lambda U'\subseteq U$ and $\lambda W'\subseteq W$. From the definition of the set ${\rm{S}}(T)$ we have that $(u+W')\cap T(U')\neq\emptyset$ and therefore $(\lambda u+W)\cap T(U)\neq \emptyset$. This means that $\lambda u\in{\rm{S}}(T)$.

The fact that ${\rm{S}}(T)$ is closed in $\mF$ is clear since it is the intersection of a family of closed subsets of $\mF$. We finally prove that if $v\in{\rm{S}}(T)$ then $(0,v)$ adheres to the graph of $T$. For any $U\in\mU_0$ and for any balanced neighborhood $V$ of $0$ in $\mF$ we find $u\in U$ and $v'\in V$ such that $v+v'=T(u)$. In other words, $-v+T(u)\in V$. Hence, $((0,v)+(U,V))\cap {\rm{Graph}}(T)\neq\emptyset$. Conversely, if $(0,v)\in\overline{{\rm{Graph}}(T)}$ then for all $U\in\mU_0$ and $V\in\mV_0$ we have $((0,v)+(U,V))\cap{{\rm{Graph}}}(T)\neq\emptyset$. It follows that $(v+V)\cap T(U)\neq\emptyset$ and then $v\in{\rm{S}}(T)$.
\end{proof}

\begin{remark}
\label{remark_ST}
Note that $(u,v)\in\overline{{\rm{Graph}}(T)}$ if and only if $(0,v-T(u))\in\overline{{\rm{Graph}}(T)}$. Indeed, we have that
\[
(u,v)\in\overline{{\rm{Graph}}(T)}\quad \Leftrightarrow\quad (0,v)\in\overline{{\rm{Graph}}(T(\cdot+u))}\quad \Leftrightarrow \quad (0,v)\in\overline{{\rm{Graph}}(T)}+(0,T(u)).
\]
It follows that the graph of $T$ is closed if and only if ${\rm{S}}(T)=\{0\}$.
\end{remark}

\begin{proposition}
\label{prop_weak_sing}
Let $(\G,\tau_1)$ be a topological $\wt{\C}$-module and $(\mF,\tau_2)$ a Hausdorff topological $\wt{\C}$-module. Every closed $\wt{\C}$-linear map $T:\G\to\mF$ has ${\rm{Ker}}(T)={\rm{K}}(T)$ and thus is weakly singular.
\end{proposition}
\begin{proof}
By Lemma \ref{lemma_pre_2}$(iii)$ we know that, since the map $T:(\G,\tau_1)\to(\mF,\tau_2)$ is closed, the strings $T(\mU)+\mV$, where $\mU$ is any topological string in $(\G,\tau_1)$ and $\mV$ is any topological string in $(\mF,\tau_2)$, generate a Hausdorff $\wt{\C}$-linear topology on $\mF$. Hence, 
\[
\bigcap_{U\in\mU_0, V\in\mV_0}T(U)+V=\{0\}.
\]
We can now write
\[
\text{Ker}(T)=T^{-1}(\cap_{U,V}(T(U)+V))=\cap_{U,V}T^{-1}(T(U)+V)=\cap_{U,V}U+T^{-1}(V)={\rm{K}}(T).
\]
It follows that ${\rm{Ker}}(T)$ is closed and coincides with ${\rm{K}}(T)$.

\end{proof}
In the sequel we denote the quotient map of $(\mF,\tau_2)$ onto $(\mF/{\rm{S}}(T),\widehat{\tau_2})$ by $\pi_T$, where $\widehat{\tau_2}$ is the quotient topology induced by $\mF$ on $\mF/{\rm{S}}(T)$. Since ${\rm{S}}(T)$ is a closed $\wt{\C}$-submodule of $\mF$, the factor $(\mF/{\rm{S}}(T),\widehat{\tau_2})$ has the structure of a Hausdorff topological $\wt{\C}$-module. 


\begin{definition}
\label{def_reg_contr}
The mapping $\pi_T\circ T:(\G,\tau_1)\to(\mF/{\rm{S}}(T),\widehat{\tau_2})$ is called \emph{the regular contraction of $T$}.
\end{definition}

\begin{proposition}
\label{prop_reg_contr}
Let $T$ be a $\wt{\C}$-linear map between the topological $\wt{\C}$-modules $(\G,\tau_1)$ and $(\mF,\tau_2)$. Then the corresponding regular contraction $\pi_T\circ T:(\G,\tau_1)\to(\mF/{\rm{S}}(T),\widehat{\tau_2})$ is closed.
\end{proposition}
\begin{proof}
By Remark \ref{remark_ST} it suffices to prove that ${\rm{S}}(\pi_T\circ T)=[0]$ in $\mF/{\rm{S}}(T)$. By definition of the singularity of $\pi_T\circ T$ we have 
\[
{\rm{S}}(\pi_T\circ T)=\cap_{U}\overline{\pi_T\circ T(U)}=\cap_U\pi_T(\overline{T(U)})=\pi_T(\cap_U\overline{T(U)})=\pi_T({\rm{S}}(T))=[0].
\]
\end{proof}

The map $T$ is related to its contraction in the following way.
\begin{proposition}
\label{prop_open_contr}
A $\wt{\C}$-linear map $T:(\G,\tau_1)\to(\mF,\tau_2)$ between topological $\wt{\C}$-modules is a strict morphism if and only if it is weakly singular and its regular contraction is a strict morphism.
\end{proposition}
\begin{proof}
By Proposition \ref{prop_strict_mor_ker} if $T$ is a strict morphism then it is weakly singular. It remains to prove that $\pi_T\circ T$ is a strict morphism, i.e., for all neighborhoods $U$ of $0$ in $(\G,\tau_1)$ the image $(\pi_T\circ T)(U)$ is a neighborhood of $0$ in $(\pi_T\circ T)(\G)$ endowed with the topology induced by $(\mF/{\rm{S}}(T),\widehat{\tau_2})$. Let us take neighborhoods $U,U'$ of $0$ in $\G$ such that $U'+U'\subseteq U$. Since $T$ is a strict morphism there exists a neighborhood $V$ of $0$ in $\mF$ such that $V\cap T(\G)\subseteq T(U')$. Finally, take a balanced neighborhood $V'$ of $0$ in $\mF$ with the property $V'+V'\subseteq V$. We want to show that $\pi_T(V')\cap \pi_T(T(\G))\subseteq \pi_T(T(U))$. The choice of neighborhoods above and the fact that $V'+{\rm{S}}(T)\subseteq V+T(U')$ yields
\begin{multline*}
\pi_T(V')\cap \pi_T(T(\G))=\pi_T((V'+{\rm{S}}(T))\cap T(\G))\subseteq \pi_T((V+T(U'))\cap T(\G))\\
=\pi_T((V\cap T(\G))+T(U'))\subseteq \pi_T(T(U)).
\end{multline*}
Conversely, assume that $T$ is weakly singular and that the regular contraction $\pi_T\circ T$ is a strict morphism. We have that ${\rm{S}}(T)\cap T(\G)=T(\cap_{V\in\mV_0}\overline{T^{-1}(V)})=T(\overline{{\rm{Ker}}(T)})=\cap_{U\in\mU_0}T(U)$. Given a neighborhood $U$ of $0$ in $\G$ we want to prove that $T(U)$ is a neighborhood of $0$ in $T(\G)$, or in other words that there exists a neighborhood $V$ of $0$ in $\mF$ such that $V\cap T(\G)\subseteq T(U)$. Let $U'$ be a neighborhood of $0$ in $\G$ such that $U'+U'\subseteq U$. Since $\pi_T\circ T$ is a strict morphism we find a neighborhood $V$ of $0$ in $\mF$ such that $\pi_T(V)\cap \pi_T(T(\G))\subseteq \pi_T(T(U'))$. It follows that $V\cap T(\G)\subseteq T(U')+({\rm{S}}(T)\cap T(\G))\subseteq T(U')+\cap_{W\in\mU_0}T(W)\subseteq T(U')+T(U')\subseteq T(U)$. This chain of inclusions completes the proof.

\end{proof}

We are now ready to work out a generalization of Theorem \ref{theorem_open_mapping_bar}. Given the topological $\wt{\C}$-modules $(\G,\tau_1)$ and $(\mF,\tau_2)$, where $\tau_1$ is a Hausdorff $\wt{\C}$-linear topology, and the surjective $\wt{\C}$-linear map $T:(\G,\tau_1)\to(\mF,\tau_2)$ we consider the regular contraction $\pi_T\circ T:(\G,\tau_1)\to (\mF/{\rm{S}}(T),\widehat{\tau_2})$ and the canonical bijection 
\[
\overline{\pi_T\circ T}:\biggl(\frac{\G}{{\rm{Ker}}(\pi_T\circ T)},\widehat{\tau_1}\biggr)\to\biggl(\frac{\mF}{{\rm{S}}(T)},\widehat{\tau_2}\biggr):u+{\rm{Ker}}(\pi_T\circ T)\to \pi_T\circ T(u)=T(u)+{\rm{S}}(T).
\]
By Proposition \ref{prop_reg_contr} the map ${\pi_T\circ T}$ is closed. Hence, by Proposition \ref{prop_weak_sing} the kernel of ${\pi_T\circ T}$  is closed and by definition of the quotient topology one can easily see that $\overline{\pi_T\circ T}$ is closed. We recall that by Proposition \ref{prop_barrelled_quo} the factor space $(\mF/{\rm{S}}(T),\widehat{\tau_2})$ is a Hausdorff barrelled $\wt{\C}$-module when $\mF$ is itself barrelled. At this point if $({\G}/{\rm{Ker}}(\pi_T\circ T),\widehat{\tau_1})$ were an infra-s $\wt{\C}$-module an application of Theorem \ref{theo_closed_barrelled_2} to $(\overline{\pi_T\circ T})^{-1}$ would allow  to conclude that $(\overline{\pi_T\circ T})^{-1}$ is continuous and therefore $\overline{\pi_T\circ T}$ is open. It would immediately follow that $\pi_T\circ T$ is open. Under the additional hypothesis that $T$ is weakly singular, Proposition \ref{prop_open_contr} would imply that $T:(\G,\tau_1)\to(\mF,\tau_2)$ is open.

This argument inspires the following definition and leads to Theorem \ref{theo_open_s}.

\begin{definition}
\label{def_s_modules}
A Hausdorff topological $\wt{\C}$-module $(\G,\tau)$ is called an \emph{s-$\wt{\C}$-module} if for every closed $\wt{\C}$-submodule $M$ the quotient $(\G/M,\widehat{\tau})$ is an infra-s $\wt{\C}$-module.
\end{definition}

\begin{theorem}
\label{theo_open_s}
Every weakly singular $\wt{\C}$-linear map from an s-$\wt{\C}$-module $(\G,\tau_1)$ onto a barrelled topological $\wt{\C}$-module $(\mF,\tau_2)$ is open.
\end{theorem}

By definition any s-$\wt{\C}$-module is an infra-s $\wt{\C}$-module. Finally, a combination of Proposition \ref{prop_quotient_com} with the considerations after Theorem \ref{theo_closed_barrelled_2} provides the following example.
\begin{proposition}
\label{prop_Fre_s}
Every Fr\'echet $\wt{\C}$-module is an s-${\wt{\C}}$-module.
\end{proposition}

As for infra-s $\wt{\C}$-modules the statement of Theorem \ref{theo_open_s} gives a characterization for s-$\wt{\C}$-modules.
\begin{proposition}
\label{prop_charac_s}
Let $(\G,\tau_1)$ be a Hausdorff topological $\wt{\C}$-module with the property that every weakly singular $\wt{\C}$-linear map from $(\G,\tau_1)$ onto a barrelled topological $\wt{\C}$-module is open. Then $(\G,\tau_1)$ is an s-$\wt{\C}$-module. 
\end{proposition}
\begin{proof}
We have to show that every quotient $(\G/M,\widehat{\tau_1})$, where $M$ is a closed $\wt{\C}$-submodule of $\G$, is an infra-s $\wt{\C}$-module. We can use the characterization of infra-s $\wt{\C}$-modules given in Remark \ref{rem_charac_infra}. Hence, it suffices to prove that every closed $\wt{\C}$-linear bijection $A$ from a Hausdorff barrelled $\wt{\C}$-module $(\mathcal{H},\tau)$ into $(\G/M,\widehat{\tau_1})$ is continuous. Let $\pi$ be the projection of $\G$ onto $\G/M$. Since $A$ is closed, $A^{-1}$ and $A^{-1}\circ\pi$ are closed maps. In particular $A^{-1}\circ\pi$ is a closed map from $(\G,\tau_1)$ into a barrelled and Hausdorff topological $\wt{\C}$-module. From Proposition \ref{prop_weak_sing} we have that $A^{-1}\circ\pi$ is weakly singular and then by hypothesis it is open. By definition of the quotient topology it follows that $A^{-1}$ is open and thus $A$ is continuous. 
\end{proof}
\begin{remark}
\label{rem_charac_s}
By the proof of the previous proposition it is clear that the following characterization holds: a Hausdorff topological $\wt{\C}$-module $(\G,\tau_1)$ is an s-$\wt{\C}$-module if and only if every closed $\wt{\C}$-linear map from $(\G,\tau_1)$ onto a barrelled and Hausdorff topological $\wt{\C}$-module is open.
\end{remark}

\section{Applications of the closed graph and the open mapping theorems}
\label{sec_application}
\subsection{Applications to Colombeau theory}
\subsubsection{Necessary condition for $\Ginf$-hypoellipticity}
The recent investigation of the $\Ginf$-regularity properties of generalized differential and pseudodifferential operators in the Colombeau context \cite{Garetto:04, GGO:03, GH:05, GH:05b, HO:03, HOP:05} has provided several sufficient conditions of $\Ginf$-hypoellipticity, i.e. technical hypotheses on the generalized symbol of the operator $P(x,D)$ which allow to conclude that $u\in\Ginf(\Om)$ when $Pu\in\Ginf(\Om)$. So far, the search for necessary condition for $\Ginf$-hypoellipticity has been a long-standing open problem. 

The closed graph theorem stated for Fr\'echet $\wt{\C}$-modules as a particular case of Theorem \ref{theorem_closed_1} for the first time enables us to find a necessary condition for $\Ginf$-hypoellipticity on the symbol of a partial differential operator with generalized constant coefficients.  

\paragraph{Preliminaries from Colombeau theory}
\leavevmode

We begin by recalling that the Colombeau algebras $\G(\Om)$ and $\Ginf(\Om)$ endowed with the corresponding sharp topologies as in \cite{Garetto:05a} Example 3.6 and 3.12 respectively have the structure of a Fr\'echet $\wt{\C}$-module. On $\G(\Om)$ we use the ultra-pseudo-seminorms $\{\mP_{K,j}\}_{K\Subset\Om,j\in \N}$ obtained through the valuation $\val_{K,j}((u_\eps)_\eps):=\sup\{b\in \R:\, \sup_{x\in K,|\alpha|\le j}|\partial^\alpha u_\eps(x)|=O(\eps^b)\, {\rm{as}}\, \eps\to 0\}$ and on $\Ginf(\Om)$ we use the family of ultra-pseudo-seminorms $\{\mP_{\Ginf(K)}\}_{K\Subset\Om}$ given by the valuations $\val_{\Ginf(K)}((u_\eps)_\eps):=\sup\{b\in\R:\, \forall\alpha\in\N^n\ \sup_{x\in K}|\partial^\alpha u_\eps(x)|=O(\eps^b)\}$. Clearly $\Ginf(\Om)$ is continuously embedded in $\G(\Om)$.

The ring of complex generalized numbers $\wt{\C}$ is a Banach $\wt{\C}$-module with respect to the ultra-pseudo-norm $|u|_\esp:=\esp^{-\val(u)}$, where $\val(u)=\val((u_\eps)_\eps):=\sup\{b\in \R:\, |u_\eps|=O(\eps^b)\, {\rm{as}}\, \eps\to 0\}$.

For more details concerning the properties of ultra-pseudo-seminorms and valuations the reader can refer to \cite{Garetto:05a}.

\begin{remark}
\label{rem_colombeau}
The spaces of Colombeau type quoted above give concrete examples of bornological, ultrabornological, barrelled and webbed $\wt{\C}$-modules as introduced in Sections \ref{section_wilde} and \ref{section_adasch}. More precisely, from Propositions 2.9, 2.14 and 2.15 in \cite{Garetto:05a} we have that $\Ginf(\Om)$, $\G(\Om)$ and $\Gc(\Om)$ are bornological and barrelled. Since they are in addition separated and complete, by Proposition \ref{prop_born_ultra}$(iv)$ they are all ultrabornological $\wt{\C}$-modules. By Propositions \ref{prop_Frechet} and \ref{prop_strict} we have that the algebras $\Ginf(\Om)$, $\G(\Om)$ and $\Gc(\Om)$ can be equipped with an absolutely convex web of type $\mC$. Finally, we recall that the duals $\LL(\Gc(\Om),\wt{\C})$ and $\LL(\G(\Om),\wt{\C})$ are endowed with the separated $\wt{\C}$-linear topology of uniform convergence on bounded subsets, i.e. with the topology obtained via the ultra-pseudo-seminorms $\mP_B(T):=\sup_{u\in B}|Tu|_\esp$ where $B$ is a bounded subset of $\Gc(\Om)$ and $\G(\Om)$ respectively. By Theorem \ref{theorem_dual} both the spaces $\LL(\Gc(\Om),\wt{\C})$ and $\LL(\G(\Om),\wt{\C})$ have an absolutely convex web of type $\mC$.
\end{remark}

\subsubsection*{Generalized points of log-type}
In stating and proving Theorem \ref{theom_nec_cond} we will consider the following set of generalized points.
\begin{definition}
\label{def_esp_point}
We say that $\xi\in\wt{\R}^n$ is of log-type if there exists a representative $(\xi_\eps)_\eps$ of $\xi\in\wt{\R}^n$ such that $|\xi_\eps|=O(\log(1/\eps))$.
\end{definition} 
The previous condition can equivalently be stated saying that the net $(\esp^{|\xi_\eps|})_\eps$ is moderate (see \cite[Lemma 2.3]{HO:03}).
Note that if $(\xi_\eps)_\eps$ and $(\xi'_\eps)_\eps$ are representatives of $\xi$ such that the nets $(\esp^{|\xi_\eps|})_\eps$ and $(\esp^{|\xi'_\eps|})_\eps$ belong both to $\EM$ then $(\esp^{|\xi_\eps|}-\esp^{|\xi'_\eps|})_\eps\in\Neg$. Indeed, from the equality
\[
\esp^{|\xi_\eps|}-\esp^{|\xi'_\eps|}=\esp^{|\xi'_\eps|}(\esp^{|\xi_\eps|-|\xi'_\eps|}-1)=\esp^{|\xi'_\eps|}(\esp^{\theta(|\xi_\eps|-|\xi'_\eps|)}(|\xi_\eps|-|\xi'_\eps|))
\]
and our assumptions it follows that for all $q\in\N$ 
\beq
\label{est_espon}
|\esp^{|\xi_\eps|}-\esp^{|\xi'_\eps|}|\le \eps^{-N+q}\esp^{||\xi_\eps|-|\xi'_\eps||}
\eeq
when $\eps$ is small enough. Since $||\xi_\eps|-|\xi'_\eps||\le |\xi_\eps-\xi'_\eps|$ we have that  $\esp^{||\xi_\eps|-|\xi'_\eps||}=O(1)$ and the estimate \eqref{est_espon} yields $(\esp^{|\xi_\eps|}-\esp^{|\xi'_\eps|})_\eps\in\Neg$. As a consequence, when $\xi$ is of log-type the generalized number
\[
\esp^{|\xi|}:=[(\esp^{|\xi_\eps|})_\eps],
\]
where $(\xi_\eps)_\eps$ is any representative of $\xi$ such that $(\esp^{|\xi_\eps|})_\eps$ is moderate, is well-defined in $\wt{\R}$. Finally, when $\zeta\in{\wt{\C}}^n$ then $\Im\zeta=(\Im\zeta_1,...,\Im\zeta_n)\in\wt{\R}^n$ and we can define the map $\wt{\C}^n\to\wt{\R}^n:\zeta\to\Im\zeta$. In particular if $\Im\zeta$ is of log-type then $\esp^{|\Im\zeta|}\in\wt{\R}$.

Before proceeding we observe that when ${\zeta}\in\wt{\C}^n$ and $\Im\zeta$ is of log-type then $$\esp^{-ix\zeta}:=[(\esp^{-ix\zeta_\eps})_\eps]$$
is a generalized function in $\G(\R^n)$. The moderateness of $(\esp^{-ix\zeta_\eps})_\eps$ is clear from $\zeta_\eps=a_\eps+i b_\eps$, $(\esp^{|b_\eps|})_\eps\in\EM$ and the inequalities
\[
\begin{split}
|\esp^{ix\zeta_\eps}|&=|\esp^{-ixa_\eps}||\esp^{xb_\eps}|\le \esp^{|x||b_\eps|},\\
|\partial^\alpha \esp^{ix\zeta_\eps}|&\le|\zeta_\eps|^{|\alpha|}|\esp^{ix\zeta_\eps}|.
\end{split}
\]
When $(|\zeta_\eps-\zeta'_\eps|)_\eps\in\Neg$ then we can write the equality
\[
\esp^{ix\zeta_\eps}-\esp^{ix\zeta'_\eps}=\esp^{+ix\zeta'_\eps}(\esp^{ix(\zeta_\eps-\zeta'_\eps)}-1)=\esp^{+ix\zeta'_\eps}\esp^{ix(\zeta_\eps-\zeta'_\eps)\theta}ix(\zeta_\eps-\zeta'_\eps).
\]
Hence arguing as above under the hypothesis that both the nets $(\esp^{|\Im\zeta_\eps|})_\eps$ and $(\esp^{|\Im\zeta'_\eps|})_\eps$ are moderate, we obtain that $(\esp^{ix\zeta_\eps}-\esp^{ix\zeta'_\eps})_\eps\in\Neg(\R^n)$.

\subsubsection*{$\Ginf$-hypoellipticity}

The regularity of a generalized function in $\G(\Om)$ with respect to $\Ginf(\Om)$ is measured by the following notions of \emph{$\Ginf$-singular support}. 
\begin{definition}
\label{def_sing_supports}
The $\Ginf$-singular support of $u\in\G(\Om)$ ($\singsupp_{\Ginf}\, u$) is the complement of the set of all points $x\in\Om$ such that the restriction of $u$ to some neighborhood $V$ of $x$ belongs to $\Ginf(V)$. 
\end{definition}
We recall that a partial differential operator with coefficients in $\wt{\C}$ maps $\Ginf(\Om)$ and $\G(\Om)$ continuously into themselves, respectively. Adopting the language already in use in \cite{HO:03} we introduce the following definition.
\begin{definition}
Let $P(D)$ be a partial differential operator with coefficients in $\wt{\C}$. $P(D)$ is said to be \emph{$\Ginf$-hypoelliptic} if for any open set $\Om\subseteq\R^n$ 
\beq
\label{Ginf_hyp}
\singsupp_{\Ginf}\, P(D)u=\singsupp_{\Ginf}\, u
\eeq
for all $u\in\G(\Om)$.
\end{definition}

\subsubsection*{Main theorem}
The following theorem is modelled on the well-known classical result valid for hypoelliptic partial differential operators with constant coefficients. For an overview on the topic we refer to \cite[Chapter1]{ChaPi:82} and \cite[Chapter 4]{Hoermander:63}.  

\begin{theorem}
\label{theom_nec_cond}
Let $P(D)$ be a partial differential operator with coefficients in $\wt{\C}$ and let $N(P)$ the set of all zeros of $P$ in $\wt{\C}^n$ with imaginary part of log-type. If $P(D)$ is $\Ginf$-hypoelliptic then
\[
\val(|\Re\zeta|)\ge 0
\]
for all $\zeta\in N(P)$.
\end{theorem}
\begin{proof} 
Let $\Om$ be an open subset of $\R^n$. Under the assumption of $\Ginf$-hy\-po\-el\-lip\-ti\-ci\-ty we know that the set ${\rm{Ker}}P:=\{u\in\G(\Om):\ P(D)u=0\}$ is contained in $\Ginf(\Om)$ and therefore in $\G(\Om)$. Since $\G(\Om)$ is a Fr\'echet $\wt{\C}$-module and $P(D)$ is a continuous map from $\G(\Om)$ into itself we conclude that ${\rm{Ker}}P$ is a Fr\'echet $\wt{\C}$-module when endowed with the topology of $\G(\Om)$.  

We can now define the map 
\[
F:{\rm{Ker}}P\to\Ginf(\Om):u\to\Delta u.
\]
An application of the closed graph theorem (Theorem \ref{theorem_closed_1}) allows to prove that $F$ is continuous. We have to verify that the graph of $F$ is sequentially closed. Let $(u,v)\in\overline{{\rm{Graph}}(F)}$, $u_n\to u$ in ${\rm{Ker}}P$ and $\Delta u_n\to v$ in
$\Ginf(\Om)$. By continuity properties we have that $\Delta u_n\to \Delta u$ in $\G(\Om)$. Since at the same time $\Delta u_n\to v$ in $\G(\Om)$ we conclude that $v=\Delta u$. Hence,  $(u,v)\in{\rm{Graph}}(F)$. 

In terms of valuations the continuity of $F$ can be expressed as follows: for all $K\Subset\Om$ there exists $L\Subset\Om$, $j\in\N$ and $c\in\R$ such that
\[
\val_{\Ginf(K)}(Fu)\ge c+\val_{L,j}(u)
\]
for all $u\in {\rm{Ker}}P$. This means that for all $u\in {\rm{Ker}}P$ and for all $m\in\N$,
\beq
\label{est_val_1}
\val_{\infty,K}(\Delta^m(Fu))\ge c+\val_{L,j}(u),
\eeq
where $\val_{\infty,K}(\Delta^m(Fu))=\val([(\sup_{x\in K}\Delta^m Fu_\eps(x))_\eps])$. 

Let $\zeta\in N(P)$. As observed above $w=\esp^{ix\zeta}\in\G(\Om)$ and since $P(D)(\esp^{ix\zeta})=P(\zeta)(\esp^{ix\zeta})$ we have that $w\in{\rm{Ker}}P$. Moreover, $Fw(x)=-\zeta^2\esp^{ix\zeta}$ for $\zeta^2=\zeta_1^2+...+\zeta_n^2$. The choice of $w$ in \eqref{est_val_1} yields
\[
\val_{\infty,K}((\zeta^2)^m\esp^{ix\zeta})\ge c+ \val_{L,j}(\esp^{ix\zeta})\ge c+ \val((1+|\zeta|)^j)+ \val(\esp^{a_1|\Im\zeta|}),
\]
where $a_1=\sup_{x\in L}|x|$. Since for $a_2=\sup_{x\in K}|x|$
\[
\val(|\zeta^2|^m\esp^{-a_2|\Im\zeta|})\ge\val_{\infty,K}((\zeta^2)^m\esp^{ix\zeta})
\]
we obtain that
\beq
\label{est_val_2}
\val(|\zeta^2|^m)\ge c+ \val((1+|\zeta|)^j)+\val(\esp^{a_1|\Im\zeta|})+\val(\esp^{a_2|\Im\zeta|})\ge c+\val((1+|\zeta|)^j)+2\val(\esp^{a|\Im\zeta|}),
\eeq
with $a=a_1+a_2$. 
The constant $c$, $a$ and $j$ in \eqref{est_val_2} do not depend on $\zeta\in N(P)$ and $m$. This leads to 
\[
\val(|\zeta^2|)\ge \frac{c}{m}+\frac{\val((1+|\zeta|)^j)}{m}+\frac{2\val(\esp^{a|\Im\zeta|})}{m}
\]
and therefore to $\val(|\zeta^2|)\ge 0$. Finally, from the equality 
\[
\label{eq_blea}
\zeta^2=\zeta_1^2+...+\zeta_n^2=|\Re\zeta|^2+2i\,\Re\zeta\,\Im\zeta-|\Im\zeta|^2.
\]
we have that $\val(||\Re\zeta|^2-|\Im\zeta|^2|)\ge\val(|\zeta^2|)\ge 0$. We can now conclude that
\[
\val(|\Re\zeta|^2)\ge\min\{\val(||\Re\zeta|^2-|\Im\zeta|^2|),\val(|\Im\zeta|^2)\}=0.
\] 
\end{proof}
\begin{remark}
\label{rem_nec_cond}
In this remark we collect some considerations concerning the necessary condition for $\Ginf$-hypoellipticity given in Theorem \ref{theom_nec_cond}. We recall that $(\omega_\eps)_\eps\in\C^{(0,1]}$ is said to be a slow scale net if for all $p\ge 0$ there exists $c_p>0$ such that $|\omega_\eps|^p\le c_p\eps^{-1}$ for all $\eps\in(0,1]$.
\begin{trivlist}
\item[(i)] The sufficient condition for $\Ginf$-hypoellipticity on a partial differential operator $P(D)=\sum_{|\alpha|\le m}c_\alpha D^\alpha$ with generalized constant coefficients presented in \cite[Definition 6.1]{GGO:03}, \cite[Definition 2.11]{Garetto:06a} and \cite[Theorem 3.2]{HOP:05} requires that a representative $(P_\eps)_\eps$ of the polynomial $P$ has the following property: there exist $l\le m$, a net $\omega_{1,\eps}\ge C\eps^{s}$ and some slow scale nets $(r_\eps)_\eps$, $(\omega_{2,\eps})_\eps$ such that 
\beq
\label{hyp_1}
|P_\eps(\xi)|\ge \omega_{1,\eps}\langle\xi\rangle^l
\eeq
and
\beq
\label{hyp_2}
|\partial^\alpha P_\eps(\xi)|\le \omega_{2,\eps}|P_\eps(\xi)|\langle\xi\rangle^{-|\alpha|}
\eeq
for all $\eps\in(0,1]$ and for all $\xi$ with $|\xi|\ge r_\eps$.

The second assertion of Theorem \ref{theom_nec_cond} allows us to claim that every partial differential operator $P(D)=\sum_{|\alpha|\le m}c_\alpha D^\alpha$ satisfying the condition above has all zeros $\zeta\in N(P)$ with $\val(|\Re\zeta|)\ge 0$. This result can be directly verified on all $\zeta=[(\zeta_\eps)_\eps]\in N(P)$ obtained by nets of zeros of $(P_\eps)_\eps$. Let $(P_\eps)_\eps$ satisfy the estimates \eqref{hyp_1} and \eqref{hyp_2} and let $(\zeta_\eps)_\eps$ be a representative of $\zeta\in N(P)$ such that $P_\eps(\zeta_\eps)=0$. Lemma 4.1.1 in \cite{Hoermander:63} shows that there exists a constant $C>0$ such that for all polynomials $p$ of degree $\le m$ we have 
\[
C^{-1}\le d(\xi)\sum_{\alpha\neq 0}|\partial^\alpha p(\xi)/p(\xi)|^{1/|\alpha|}\le C,\quad \xi\in\R^n,\ p(\xi)\neq 0,
\]
where $d(\xi)$ denotes the distance from $\xi$ to the surface $\{\zeta\in\C^n:\, p(\zeta)=0\}$. It follows that if $|\Re\zeta_\eps|\ge r_\eps$ then by \eqref{hyp_2} there exists some slow scale net $(\omega_\eps)_\eps$ such that
\[
C^{-1}\le d(\Re\zeta_\eps)\sum_{\alpha\neq 0}|\partial^\alpha P_\eps(\Re\zeta_\eps)/P_\eps(\Re\zeta_\eps)|^{1/|\alpha|}\le \omega_\eps d(\Re\zeta_\eps) \langle\Re\zeta_\eps\rangle^{-1}\le \omega_\eps d(\Re\zeta_\eps) |\Re\zeta_\eps|^{-1}.
\]
Hence,
\[
|\Re\zeta_\eps|\le C\omega_\eps d(\Re\zeta_\eps)\le C\omega_\eps |\Im\zeta_\eps|\le C'\omega_\eps\log(1/\eps).
\]
In conclusion $|\Re\zeta_\eps|\le \max\{r_\eps,C'\omega_\eps\log(1/\eps)\}$ for all $\eps\in(0,1]$ and therefore $\val(|\Re\zeta|)\ge 0$. 
\item[(ii)] If the polynomial $P$ has a zero $\zeta$ with classical imaginary part such that $\val(|\Re\zeta|)<0$ then the corresponding operator is not $\Ginf$-hypoelliptic. As an example take the operator $P(D)=-i[(\eps^r)_\eps]\partial_{x_1}+\partial_{x_2}$ with $r>0$. The point $\zeta=(\zeta_1,\zeta_2)$ with $\zeta_1=[(\eps^{-r}+i)_\eps]$ and $\zeta_2=[(-\eps^{r}+i)_\eps]$ has $\val(|\Re\zeta|)=-r$ and the operator $P(D)$ is not $\Ginf$-hypoelliptic. Indeed, the generalized function $u=[(\esp^{(ix_1\eps^{-r}-x_2)})_\eps]$ satisfies $P(D)u=0$ but $u\not\in\Ginf(\R^2)$.
 
\end{trivlist}
\end{remark}

\subsubsection{Cauchy problem in the Colombeau framework: continuous dependence on the data}
Let us consider the initial value problem
\beq
\label{cauchy_prob}
\dot{u}(t)=a(t)u(t)+b(t),\qquad\qquad\qquad u(t_0)=0
\eeq
in the Colombeau context, where $a,b\in\G(\R)$ and look for a solution $u$ in $\G(\R)$. By Theorem 1.5.2 and Remark 1.5.3 in \cite{GKOS:01} we know that if $a$ is of $L^\infty$-log type, that is it has a representative $(a_\eps)_\eps$ with the property $\Vert a_\eps\Vert_{L^\infty}=O(\log(1/\eps))$, then the problem \eqref{cauchy_prob} has a unique solution in $\G(\R)$. 

We define $D_{t_0}:=\{u\in\G(\R):\, u(t_0)=0\}$ and the $\wt{\C}$-linear map 
\[
T:D_{t_0}\to\G(\R):u\to \dot{u}-au.
\]
Since the evaluation map $\G(\R)\to\wt{\C}:u\to u(t_0)$ is continuous $D_{t_0}$ is a closed $\wt{\C}$-submodule of $\G(\R)$ and therefore it is a Fr\'echet $\wt{\C}$-module with the topology induced by $\G(\R)$. One can easily see that the $\wt{\C}$-linear bijective map $T$ is continuous. 
We can then apply the open mapping theorem \ref{theorem_open} (which is actually in the form of an isomorphism theorem in this case) and conclude that the map $T$ has continuous inverse. This means that the solution $u$ of the Cauchy problem \eqref{cauchy_prob} depends continuously on $b$.

More generally the continuous dependence on the inhomogeneity $b$ and the initial values $u_0, u_1,...,u_{n-1}\in\wt{\R}$ holds for any Cauchy problem uniquely solvable in $\G(\R)$ of the form 
\[
u^{(n)}(t)=\sum_{i=1}^{n-1}a_i(t)u^{(i)}(t)+b(t),\qquad\quad u(\wt{t_0})=u_0,\ u'(\wt{t_0})=u_1,\,...,\,u^{(n-1)}(\wt{t_0})=u_{n-1}
\]
where $a_i,b\in\G(\R)$, $i=1,...,n-1$ and $\wt{t_0}$ is a generalized point with compact support. This is obtained by proving that the inverse of the bijection 
\[
T:\G(\R)\to\G(\R)\times\wt{\R}^n:u\to \biggl( u^{(n)}-\sum_{i=1}^{n-1}a_iu^{(i)},u(\wt{t_0}), u'(\wt{t_0}),...,u^{(n-1)}(\wt{t_0})\biggr)
\]
is continuous via application of the open mapping theorem for $\wt{\C}$-linear maps between Fr\'echet $\wt{\C}$-modules.

This functional analytic method is particularly convenient when we deal with Cauchy problems given by a partial differential equation because the direct investigation of the continuous dependence on the data of the solution may involve rather complicated estimates. As an instructive example let us consider the linear wave equation 
\beq
\label{cauchy_wave}
(\partial^{2}_{t}-\Delta)u=f\qquad \text{in}\, \G(\R^{n+1}),\qquad u|_{\{t=0\}}=c,\qquad \partial_{t} u|_{\{t=0\}}=d\qquad  \text{in}\, \G(\R^{n}), 
\eeq
which has a unique solution in $\G(\R^{n+1})$ given $c,d\in\G(\R^n)$ and $f\in\G(\R^{n+1})$ (see \cite[Remark 16.4]{O:92}). The map 
\[
T:\G(\R^{n+1})\to\G(\R^n)\times\G(\R^n)\times\G(\R^{n+1}):u\to (u|_{\{t=0\}},\partial_{t} u|_{\{t=0\}},(\partial^{2}_{t}-\Delta)u)
\]
is bijective and continuous. By application of the open mapping theorem \ref{theorem_open} we conclude that the solution $u\in\G(\R^{n+1})$ of the Cauchy problem \eqref{cauchy_wave} depends continuously on $f$, $c$ and $d$.   

\subsection{Applications to the theory of Banach $\wt{\C}$-modules}

We conclude the paper by presenting some applications of the closed graph and the open mapping theorems elaborated in Sections \ref{section_wilde}, \ref{section_open}, \ref{section_adasch} to the abstract theory of Banach $\wt{\C}$-modules (i.e., complete and ultra-pseudo-normed $\wt{\C}$-modules).

\subsubsection{Equivalence of ultra-pseudo-norms}
\begin{definition}
\label{def_eq_norm}
The ultra-pseudo-norms $\mP$ and $\mQ$ on the $\wt{\C}$-module $\G$ are said to be \emph{equivalent} if there exist some constants $C_1,C_2>0$ such that $C_1\mP(u)\le\mQ(u)\le C_2\mP(u)$ for all $u\in\G$.
\end{definition}

\begin{proposition}
\label{prop_eq_norm}
Let $\G$ be a $\wt{\C}$-module and $\mP$ and $\mQ$ two ultra-pseudo-norms which give the structure of a Banach $\wt{\C}$-module to $\G$. If $\mQ(u)\le C\mP(u)$ for some $C>0$ and for all $u\in\G$ then $\mP$ and $\mQ$ are equivalent.
\end{proposition}
\begin{proof}
By the inequality involving $\mP$ and $\mQ$ we know that the identity map from $(\G,\mP)$ into $(\G,\mQ)$ is continuous. An application of Theorem \ref{theorem_open} yields the continuity of the identity map from $(\G,\mQ)$ to $(\G,\mP)$. As a consequence, $\mP$ and $\mQ$ are equivalent ultra-pseudo-norms on $\G$.
\end{proof}

\subsubsection{A theorem of uniform boundedness}
The Banach $\wt{\C}$-module-version of the closed graph theorem allows to obtain a result of uniform boundedness for $\wt{\C}$-linear maps. We begin by recalling some useful preliminary notions of functional analysis.
\begin{proposition}
\label{prop_pre_fa}
Let $(\G,\mP)$ and $(\mF,\mQ)$ be ultra-pseudo-normed $\wt{\C}$-modules.
\begin{itemize}
\item[(i)] The $\wt{\C}$ module
\[
\LL(\G,\mF):=\{T:\G\to\mF,\ \wt{\C}-\text{linear and continuous}\}
\]
can be equipped with the structure of an ultra-pseudo-normed $\wt{\C}$-module by means of the ultra-pseudo-norm \[
\mP_{\LL(\G,\mF)}(T):=\sup_{\mP(u)=1}\mQ(Tu).
\]
\item[(ii)] If $(\mF,\mQ)$ is complete then $(\LL(\G,\mF),\mP_{\LL(\G,\mF)})$ is complete.
\item[(iii)] Let $X$ be a subset of $\G$. Then the set
\[
\mB(X,\mF):=\{T:X\to\mF,\ \sup_{u\in X}\mQ(Tu)<\infty\}
\]
is an ultra-pseudo-normed $\wt{\C}$-module with respect to the ultra-pseudo-norm
\[
\mP_{\mB(X,\mF)}(T):=\sup_{u\in X}\mQ(Tu).
\]
\item[(iv)] If $(\mF,\mQ)$ is complete then $(\mB(X,\mF),\mP_{\mB(X,\mF)})$ is complete.
\end{itemize}
\end{proposition}
\begin{proof}
We leave it to the reader to prove the assertions $(i)$ and $(ii)$ since they are analogous to the classical statements valid for normed and Banach spaces. The same holds for $(iii)$. Assume now that $(\mF,\mQ)$ is complete. We want to prove that any Cauchy sequence $(T_n)_n$ in $(\mB(X,\mF),\mP_{\mB(X,\mF)})$ is convergent. Since for all $\delta>0$ there exists $N\in\N$ such that $\mP_{\mB(X,\mF)}(T_n-T_m)\le\delta$ for all $n,m\ge N$ it follows that $(T_n(u))_n$ is a Cauchy sequence in $\mF$ for all fixed $u\in X$. Consequently it converges to some $Tu$ in $\mF$. Combining the assertion
\[
\forall u\in X\ \forall \delta>0\ \exists N\in\N\ \forall n\ge N\qquad \mQ(Tu-T_n u)\le\delta
\]
with
\[
\forall \delta>0\  \exists M\in\N\ \forall m,n\ge M\ \forall u\in X\qquad \mQ(T_n u-T_mu)\le \delta
\]
we obtain that $\sup_{u\in X}\mQ(Tu-T_nu)\le\delta$ for $n$ large enough. It follows that $T\in\mB(X,\mF)$ and that $T_n\to T$ in $\mB(X,\mF)$.
\end{proof}

\begin{theorem}
\label{theo_uni_bound}
Let $(\G,\mP)$ and $(\mF,\mQ)$ be Banach $\wt{\C}$-modules. Let $Y\subseteq \LL(\G,\mF)$ be a pointwise bounded subset of $\LL(\G,\mF)$, i.e., $\sup_{T\in Y}\mQ(Tu)<\infty $ for all $u\in\G$. Then $Y$ is a bounded subset of the $\wt{\C}$-module $(\LL(\G,\mF),\mP_{\LL(\G,\mF)})$.
\end{theorem}
\begin{proof}
By Proposition \ref{prop_pre_fa} the $\wt{\C}$-module $\mB(Y,\mF)$ endowed with the topology induced by the ultra-pseudo-norm $\mP_{\mB(Y,\mF)}$ is a Banach $\wt{\C}$-module. We define the map
\[
S:(\G,\mP)\to (\mB(Y,\mF),\mP_{\mB(Y,\mF)}):u\to (T\to T(u)).
\]
This is possible since by the assumption of pointwise boundedness of $Y$ we have $\sup_{T\in Y}\mQ(S(u)(T))<\infty$. The map $S$ is $\wt{\C}$-linear and its graph is sequentially closed. Indeed, if $u_n\to u$ in $(\G,\mP)$ and $Su_n\to H$ in $\mB(Y,\mF)$ then for all $T\in Y$ we have $\mQ(S(u_n)(T)-H(T))\to 0$. This means that $Tu_n\to H(T)$ in $\mF$. But by the continuity of $T$ we also have that $Tu_n\to Tu$ in $\mF$. Hence $H(T)=T(u)$ for all $T\in Y$ or in other words $H=S(u)$. We can now apply Theorem \ref{theorem_closed_1} and conclude that the map $S$ is continuous. It follows that there exists a constant $C>0$ such that $\sup_{\mP(u)=1}\mP_{\mB(Y,\mF)}(S(u))\le C$ and therefore for all $T\in Y$,
\[
\sup_{\mP(u)=1}\mQ(Tu)=\sup_{\mP(u)=1}\mQ(S(u)(T))\le\sup_{\mP(u)=1}\sup_{T\in Y}\mQ(S(u)(T))=\sup_{\mP(u)=1}\mP_{\mB(Y,\mF)}(S(u))\le C.
\]
In conclusion, $\sup_{T\in Y}\mP_{\LL(\G,\mF)}(T)\le C$.

\end{proof}

\begin{remark}
\label{rem_Ban_bar}
More generally one can prove that if $(\G,\tau_1)$ and $(\mF,\tau_2)$ are topological $\wt{\C}$-modules then a pointwise bounded subset $Y$ of $\LL(\G,\mF)$ is equicontinuous with respect to $(\G,\tau_1^b)$ and $(\mF,\tau_2)$. Indeed, if $\mV$ be a closed topological string in $(\mF,\tau_2)$ then $T^{-1}(\mV)$ is a closed topological string in $(\G,\tau_1)$ for any $T\in Y$. From the pointwise boundedness of $Y$ it follows that $\mU:=\cap_{T\in Y}T^{-1}(\mV)$ is a closed topological string in $(\G,\tau_1)$. Hence, $\mU$ is a topological string in $(\G,\tau_1^b)$. Finally, since by construction $T(\mU)\subseteq \mV$ for every $T\in Y$, we conclude that $Y$ is equicontinuous with respect to $(\G,\tau_1^b)$ and $(\mF,\tau_2)$.\\
In particular, $Y$ is an equicontinuous subset of $\LL((\G,\tau_1),(\mF,\tau_2))$ when the $\wt{\C}$-module $(\G,\tau_1)$ is barrelled.

\end{remark}
\section*{Appendix: bornological and ultrabornological $\wt{\C}$-modules} 

\setcounter{section}{1}
\renewcommand{\thesection}{\Alph{section}}
\renewcommand{\theequation}{A.\arabic{equation}}
\setcounter{equation}{0}

In this appendix we give a brief survey on bornological and ultrabornological $\wt{\C}$-module. We begin with the definition of a bornological $\wt{\C}$-module given in \cite[Definition 2.7]{Garetto:05a}.
\begin{definition}
\label{def_born_app}
A subset $A$ of a locally convex topological $\wt{\C}$-module $\G$ is said to be \emph{bornivorous} if it absorbs any bounded subset of $\G$.
A locally convex topological $\wt{\C}$-module $\G$ is \emph{bornological} if every absolutely convex and bornivorous subset of $\G$ is a neighborhood of the origin. 
\end{definition}
Every ultra-pseudo-normed $\wt{\C}$-module is bornological or more generally, as proved in \cite[Proposition 2.9]{Garetto:05a}, every locally convex topological $\wt{\C}$-module which has a countable base of neighborhoods of the origin is bornological.


It is possible to characterize the topology of a  Hausdorff bornological $\wt{\C}$-module as the finest locally convex $\wt{\C}$-linear topology which makes any injection $\iota_B:\G_B\to \G$ continuous, where $B$ is any bounded disk. An analogous conclusion holds for Hausdorff ultrabornological $\wt{\C}$-modules when $B$ is a bounded Banach disk. To achieve these characterizations a deeper investigation of the structural properties of the $\wt{\C}$-module $\G_A$ is needed.

\begin{proposition}
\label{prop_G_A}
Let $\G$ be a locally convex topological $\wt{\C}$-module and $A$ a nonempty, balanced and convex subset of $\G$.
\begin{itemize}
\item[(i)] If $\G$ is a Hausdorff $\wt{\C}$-module and $A$ is bounded then $\G_A$ is an ultra-pseudo-normed $\wt{\C}$-module.
\item[(ii)] If in addition to the hypotheses of $(i)$ the subset $A$ is complete, then $\G_A$ is a Banach $\wt{\C}$-module.
\end{itemize}
\end{proposition}
 
The proof of the second assertion of Proposition \ref{prop_G_A} makes use of the following lemma.
\begin{lemma}
\label{lemma_complete}
Let $\G$ be a $\wt{\C}$-module and let $\tau$ and $\tau'$ be two Hausdorff $\wt{\C}$-linear topologies on $\G$. Suppose that $\tau$ is finer than $\tau'$ and that it has a fundamental system $\mathcal{B}$ of balanced neighborhoods of the origin which are complete for $\tau'$. Then $\G$ is complete for $\tau$.
\end{lemma}
\begin{proof}
We begin by proving that every $V\in\mathcal{B}$ is complete for $\tau$. Let $\mF$ be a Cauchy filter on $V$. Since $\mF$ is also a Cauchy filter with respect to the topology $\tau'$, we find some $x_0\in V$ such that $\mF\to x_0$ for $\tau'$ on $V$. Let $W$ be a neighborhood of $0$ for $\tau$ such that $W$ is balanced and complete for $\tau'$ and $W+W\subseteq V$. It follows that $W$ is closed for $\tau'$. By definition of a Cauchy filter there exists $A\in \mF$ such that $A-A\subseteq W$. Let us take $x_1\in A$. Then $A\subseteq x_1+W$ and since $x_1+W$ is closed for $\tau'$, the closure of $A$ with respect to $\tau'$ is contained in $x_1+W$. Therefore, $x_0\in x_1+W$ and by the fact that $W$ is balanced we conclude that $A\subseteq x_0+W+W\subseteq x_0+V$. This proves that $\mF\to x_0$ for $\tau$. Finally, assume that $\mF$ is a Cauchy filter on $\G$ for $\tau$. For all $V\in\mathcal{B}$ there exists $A\in\mF$ such that $A-A\subseteq V$. Let us fix $V$ and $A$ with the previous property. For $x_1\in A$ we have that $A\subseteq x_1+V$ and therefore the filter induced by $\mF$ on $x_1+V$ is a Cauchy filter. Since $V$ is complete with respect to the topology $\tau$, $x_1+V$ is also complete and therefore the filter induced by $\mF$ on $x_1+V$ is convergent. It follows that $\mF$ is convergent for $\tau$ and $\G$ is complete. 
\end{proof}
\begin{proof}[Proof of Proposition \ref{prop_G_A}]
$(i)$ Let $u\neq 0$ be an element of $\G_A$. We want to prove that $\mP_A(u)\neq 0$. Since $\G$ is Hausdorff there exists a neighborhood $U$ of the origin such that $u\not\in U$ and by the boundedness of $A$ we know that there exists $a\in\R$ such that $[(\eps^{-b})_\eps]A\subseteq U$ for all $b\le a$. Consequently $u\not\in[(\eps^{-b})_\eps]A$ and $\val_A(u)\neq +\infty$.

$(ii)$ In the sequel we denote the topology generated by the ultra-pseudo-norm $\mP_A$ on $\G_A$ with $\tau$ and the topology induced by $\G$ on $\G_A$ with $\tau'$. Since for all $\eta>0$,
\[
\{u\in \G_A:\, \mP_A(u)<\eta\}\subseteq [(\eps^{-\log\eta})_\eps]A\subseteq\{u\in \G_A:\, \mP_A(u)\le\eta\},
\]
$([(\eps^n)_\eps]A)_n$ is a fundamental system of balanced neighborhoods of the origin for $\tau$. Each of them is complete for $\tau'$ by the completeness of $A$. Since $A$ is bounded, for every neighborhood $U$ of $0$ in $\G$ there exists $a$ such that $[(\eps^{-b})_\eps]A\subseteq U$ for all $b\le a$. This means that $\tau$ is finer than $\tau'$. Therefore, by applying Lemma \ref{lemma_complete} we conclude that $\G_A$ is complete for the topology $\tau$. 
\end{proof}

A straightforward consequence of Proposition \ref{prop_G_A} is the following result on bornological $\wt{\C}$-modules.  \begin{proposition}
\label{prop_born_1}
Let $\G$ be a bornological locally convex topological $\wt{\C}$-module. Then there exists a family of locally convex topological $\wt{\C}$-modules $(\G_\gamma)_{\gamma\in \Gamma}$, whose topology can be defined with the help of a single ultra-pseudo-seminorm, and a family $(\iota_\gamma)_{\gamma\in \Gamma}$ of $\wt{\C}$-linear maps $\iota_\gamma:\G_\gamma\to \G$ such that $\G$ is endowed with the finest locally convex $\wt{\C}$-linear topology which makes every $\iota_\gamma$ continuous. Moreover, if $\G$ is a Hausdorff $\wt{\C}$-module then every $\G_\gamma$ can be taken to be an ultra-pseudo-normed $\wt{\C}$-module. Finally, if $\G$ is Hausdorff and quasi-complete then every $\G_\gamma$ can be taken to be a Banach $\wt{\C}$-module.
\end{proposition}
\begin{proof}
Let $A$ be an absolutely convex closed and bounded subset of $\G$. We know that $\G_A$ is topologized through the ultra-pseudo-seminorm $\mP_A$. Let $\iota_A$ be the canonical injection of $\G_A$ into $\G$.
We denote the original topology on $\G$ by $\tau$ and the finest locally convex $\wt{\C}$-linear topology which makes every $\iota_A$ continuous by $\tau'$. Let $V$ be an absolutely convex neighborhood of $0$ for $\tau$. Then there exists $a\in\R$ such that $[(\eps^{-a})_\eps]A\subseteq V$. This means that $[(\eps^{-a})_\eps]A\subseteq \iota_{A}^{-1}(V)$ and therefore $V$ is a neighborhood of the origin in $\G$ with respect to $\tau'$. Conversely, let $V$ be an absolutely convex neighborhood of the origin for $\tau'$. Then for each absolutely convex bounded and closed subset $A$ of $\G$, $\iota_{A}^{-1}(V)$ is a neighborhood of $0$ for $\G_A$. Hence for every $A$ there exists $a\in\R$ such that $[(\eps^{-a})_\eps]A\subseteq \iota_{A}^{-1}(V)$ and then $[(\eps^{-a})_\eps]A\subseteq V$. It follows that $V$ is an absolutely convex and bornivorous subset of $\G$. Since $\G$ is bornological we conclude that $V$ is a neighborhood of $0$ for $\tau$. The other two assertions of this proposition are entailed by Proposition \ref{prop_G_A}.
\end{proof}


\begin{proposition}
\label{prop_born_ultra}
\leavevmode
\begin{itemize}
\item[(i)] A Hausdorff locally convex topological $\wt{\C}$-module $\G$ is bornological if and only if it has the finest locally convex $\wt{\C}$-linear topology given by the family of ultra-pseudo-normed $\wt{\C}$-modules $(\G_B)_B$, where $B$ is any bounded disk.
\item[(ii)] Every ultrabornological $\wt{\C}$-module $\G$ is bornological.
\item[(iii)] A Hausdorff locally convex topological $\wt{\C}$-module $\G$ is ultrabornological if and only if it has the finest locally convex $\wt{\C}$-linear topology determined by the family of Banach $\wt{\C}$-modules $(\G_B)_B$, where $B$ is any bounded Banach disk if and only if it has the finest locally convex $\wt{\C}$-linear topology determinated by a family of $\wt{\C}$-linear injections $\iota_\gamma:\G_\gamma\to\G$, where every $\G_\gamma$ is a Banach $\wt{\C}$-module.
\item[(iv)] Every bornological, Hausdorff, quasi-complete $\wt{\C}$-module is ultrabornological.
\end{itemize}
\end{proposition}
\begin{proof}
$(i)$ By Proposition \ref{prop_born_1} the necessary condition is clear. Conversely assume that $\G$ has the finest locally convex $\wt{\C}$-linear topology which makes every map $\iota_B:\G_B\to\G$ continuous. We have to prove that every absolutely convex and bornivorous subset $V$ of $\G$ is a neighborhood of the origin. Clearly, for every $B$ we have that $\iota_B^{-1}(V)\subseteq\G_B$ is absolutely convex and since by assumption $\iota_B$ is continuous, every bounded subset $D$ of $\G_B$ is mapped into a bounded subset of $\G$. $V$ is bornivorous. This means that given $D\subseteq \G_B$ there exists $d'\in\R$ such that $\iota_B(D)\subseteq[(\eps^d)_\eps]V$ for all $d'\le d$, that is $D\subseteq[(\eps^d)_\eps]\iota_B^{-1}(V)$ for all $d'\le d$. Hence, by the fact that $\G_B$ is ultra-pseudo-normed and therefore bornological, we conclude that $\iota_B^{-1}(V)$ is a neighborhood of the origin in $\G_B$. This proves that $V$ is a neighborhood of the origin in $\G$ and that $\G$ is bornological.

$(ii)$ We have to prove that if $\G$ is ultrabornological then every absolutely convex and bornivorous subset $U$ of $\G$ is a neighborhood of the origin. Clearly, since $U$ is bornivorous it is also ultrabornivorous and therefore it is a neighborhood of $0$ in $\G$. 

 
$(iii)$ As in Proposition \ref{prop_born_1} one easily proves that if $\G$ is ultrabornological then it has the finest locally convex $\wt{\C}$-linear topology determined by the family of injections $\iota_B:\G_B\to \G$ where $B$ is any bounded Banach disk. Conversely, assume that $\G$ is Hausdorff and has the finest locally convex $\wt{\C}$-linear topology which makes any map $\iota_B:\G_B\to\G$ continuous where $B$ is a bounded Banach disk. We want to prove that every absolutely convex and ultrabornivorous subset $V$ of $\G$ is a neighborhood of $0$. Let $D$ be a bounded Banach disk in $\G_B$. This means that $D$ is absolutely convex, bounded in $\G_B$ and that $((\G_B)_D,\mP_D)=(\G_D,\mP_D)$ is Banach. Hence, $\iota_B(D)$ is a Banach disk of $\G$ and since $V$ is ultrabornivorous there exists $d'\in\R$ such that $\iota_B(D)\subseteq[(\eps^d)_\eps]V$ for all $d'\le d$, that is $D\subseteq[(\eps^d)_\eps]\iota_B^{-1}(V)$ for all $d'\le d$. This shows that $\iota_B^{-1}(V)$ is an absolutely convex ultrabornivorous subset of $\G_B$. Since $\G_B$ is a Banach and therefore ultrabornological $\wt{\C}$-module we obtain that $\iota_B^{-1}(V)$ is a neighborhood of $0$ in $\G_B$. Hence, $V$ is a neighborhood of $0$ in $\G$. 

It remains to prove that if $\G$ has the finest locally convex $\wt{\C}$-linear topology determinated by the family of $\wt{\C}$-linear injections $\iota_\gamma:\G_\gamma\to\G$, where every $\G_\gamma$ is a Banach $\wt{\C}$-module, then every absolutely convex and ultrabornivorous subset $U$ of $\G$ is a neighborhood of the origin. By definition $\iota_\gamma^{-1}(U)$ is absolutely convex in $\G_\gamma$. Let $B$ be a bounded Banach disk in $\G_\gamma$. Then $\iota_\gamma(B)$ is a bounded disk in $\G$. We want to prove that $(\G_{\iota_\gamma(B)},\mP_{\iota_\gamma(B)})$ is a Banach $\wt{\C}$-module, knowing that $((\G_\gamma)_B,\mP_B)$ is a Banach $\wt{\C}$-submodule of $\G_\gamma$. Let $(u_n)_n$ be a Cauchy sequence in $\G_{\iota_\gamma(B)}$. Then for all $n$ there exist $m_n\in\N$ and $b_n\in B$ such that $u_n=[(\eps^{-m_n})_\eps]\iota_\gamma(b_n)$ and $\val_{\iota_\gamma(B)}([(\eps^{-m_{n+1}})_\eps]\iota_\gamma(b_{n+1})-[(\eps^{-m_n})_\eps]\iota_\gamma(b_n))\to+\infty$. It follows that the sequence $([(\eps^{-m_n})_\eps]b_n)_n$ belongs to $(\G_\gamma)_B$ and that for all $q\in\N$ there exists $N\in\N$  such that $\iota_\gamma([(\eps^{-m_{n+1}})_\eps]b_{n+1}-[(\eps^{-m_n})_\eps]b_n)\in[(\eps^q)_\eps]\iota_\gamma(B)$ for all $n\ge N$. Thus, $[(\eps^{-m_{n+1}})_\eps]b_{n+1}-[(\eps^{-m_n})_\eps]b_n \in [(\eps^q)_\eps]B$ for all $n\ge N$ or in other words $([(\eps^{-m_n})_\eps]b_n)_n$ is a Cauchy sequence in $(\G_\gamma)_B$. Since $(\G_\gamma)_B$ is complete, $([(\eps^{-m_n})_\eps]b_n)_n$ tends to some $v\in(\G_\gamma)_B$ and by continuity of $\iota_\gamma$ we conclude that $u_n\to \iota_\gamma(v)$ in $\G_{\iota_\gamma(B)}$. This shows that $(\G_{\iota_\gamma(B)},\mP_{\iota_\gamma(B)})$ is complete. Now, $\iota_\gamma(B)$ is a bounded Banach disk in $\G$ and hence it is absorbed by $U$, i.e., $\iota_\gamma(B)\subseteq [(\eps^b)_\eps]U$ for all $b$ smaller than a certain $a\in\R$. It follows that $B\subseteq[(\eps^b)_\eps]\iota_\gamma^{-1}(U)$. We conclude that $\iota_\gamma^{-1}(U)$ is an absolutely convex ultrabornivorous subset of $\G_\gamma$ and therefore a neighborhood of $0$ in $\G_\gamma$. As a straightforward consequence we have that $U$ is a neighborhood of $0$ in $\G$.


$(iv)$ From the first assertion of this proposition combined with Proposition \ref{prop_born_1} we know that if $\G$ is bornological and Hausdorff then it has the finest locally convex $\wt{\C}$-linear topology which makes any injection $\iota_B:\G_B\to\G$ continuous, where $B$ is a closed and bounded disk of $\G$. By Proposition \ref{prop_G_A} and the assumption of quasi-completeness we have that every $(\G_B,\mP_B)$ is a Banach $\wt{\C}$-module. The third claim of this proposition allows to conclude that $\G$ is ultrabornological.
\end{proof}
Since every Fr\'echet $\wt{\C}$-module is bornological, Hausdorff and complete by the fourth statement of the previous proposition we obtain that every Fr\'echet $\wt{\C}$-module is ultrabornological.

\paragraph{Acknowledgment:} The author is grateful to Prof. Michael Oberguggenberger for several
inspiring discussions during the preparation of the paper. 
\bibliographystyle{abbrv}

\begin{thebibliography}{10}

\bibitem{Adasch:78}
N.~Adasch, B.~Ernst, and D.~Keim.
\newblock {\em Topological Vector Spaces}.
\newblock Number 639 in Lecture Notes in Math. Springer-Verlag, Berlin, 1978.

\bibitem{ChaPi:82}
J.~Chazarain and A.~Piriou.
\newblock {\em Introduction to the theory of linear partial differential
  equations}.
\newblock Number~14 in Studies in Mathematics and its Applications. North
  Holland Publishing and Co., Amsterdam, 1982.

\bibitem{DeWilde:67}
M.~{De Wilde}.
\newblock Sur le th\'eor\`em du graphe ferm\'e.
\newblock {\em C. R. Acad. Sc. Paris}, 265, s\'erie A:376--379, 1967.

\bibitem{DeWilde:71}
M.~{De Wilde}.
\newblock Ultrabornological spaces and the closed graph theorem.
\newblock {\em Bull. Soc. R. Sc. Li\`ege}, 40:116--118, 1971.

\bibitem{DeWilde:78}
M.~{De Wilde}.
\newblock {\em Closed graph theorems and webbed spaces}, volume~19.
\newblock Research Notes in Mathematics, Pitman, London, 1978.

\bibitem{Garetto:04}
C.~Garetto.
\newblock Pseudo-differential operators in algebras of generalized functions
  and global hypoellipticity.
\newblock {\em Acta Appl. Math.}, 80(2):123--174, 2004.

\bibitem{Garetto:05a}
C.~Garetto.
\newblock Topological structures in {C}olombeau algebras: topological
  $\wt{\C}$-modules and duality theory.
\newblock {\em Acta. Appl. Math.}, 88(1):81--123, 2005.

\bibitem{Garetto:05b}
C.~Garetto.
\newblock Topological structures in {C}olombeau algebras: investigation of the
  duals of ${\Gc(\Om)}$, ${\G(\Om)}$ and ${\GS(\R^n)}$.
\newblock {\em Monatsh. Math.}, 146(3):203--226, 2005.
 
\bibitem{Garetto:06a}
C.~Garetto.
\newblock Microlocal analysis in the dual of a {C}olombeau algebra: generalized
  wave front sets and noncharacteristic regularity.
\newblock {\em arXiv:math. AP/0511297}, 2006.

\bibitem{GGO:03}
C.~Garetto, T.~Gramchev, and M.~Oberguggenberger.
\newblock Pseudodifferential operators with generalized symbols and regularity
  theory.
\newblock {\em Electron. J. Diff. Eqns.}, 2005(2005)(116):1--43, 2003.

\bibitem{GH:05}
C.~Garetto and G.~H\"{o}rmann.
\newblock Microlocal analysis of generalized functions: pseudodifferential
  techniques and propagation of singularities.
\newblock {\em Proc. Edinburgh. Math. Soc.}, 48(3):603--629, 2005.

\bibitem{GH:05b}
C.~Garetto and G.~H\"ormann.
\newblock Duality theory and pseudodifferential techniques for Colombeau
  algebras: generalized kernels and microlocal analysis.
\newblock {\em Proceedings of the Conference ``Generalized Functions 2004'',
  University of Novi Sad}, 2005.
\newblock to appear in Bull. Cl. Sci. Math. Nat. Sci. Math.
 
\bibitem{GKOS:01}
M.~Grosser, M.~Kunzinger, M.~Oberguggenberger, and R.~Steinbauer.
\newblock {\em Geometric theory of generalized functions}, volume 537 of {\em
  Mathematics and its Applications}.
\newblock Kluwer, Dordrecht, 2001.

\bibitem{Hoermander:63}
L.~H{\"o}rmander.
\newblock {\em Linear Partial Differential Operators}.
\newblock Springer-Verlag, Berlin, 1963.

\bibitem{HO:03}
G.~H{\"o}rmann and M.~Oberguggenberger.
\newblock Elliptic regularity and solvability for partial differential
  equations with {C}olombeau coefficients.
\newblock {\em Electron. J. Diff. Eqns.}, 2004(14):1--30, 2004.

\bibitem{HOP:05}
G.~H\"{o}rmann, M.~Oberguggenberger, and S.~Pilipovic.
\newblock Microlocal hypoellipticity of linear partial differential operators
  with generalized functions as coefficients.
\newblock {\em Trans. Amer. Math. Soc.}, 2005.
\newblock to appear.

\bibitem{Horvath:66}
J.~Horv\'{a}th.
\newblock {\em Topological vector spaces and distributions}.
\newblock Addison-Wesley, Reading, MA, 1966.

\bibitem{Koethe:79}
G.~K\"othe.
\newblock {\em Topological vector spaces II}.
\newblock Springer, New York, 1979.

\bibitem{O:82}
M.~Oberguggenberger.
\newblock {\em Der Graphensatz in lokalkonvexen topologischen Vectorr\"aumen}.
\newblock Teubner-texte zur Mathematik. Teubner Verlagsgesellschaft, Leipzig,
  1982.

\bibitem{O:92}
M.~Oberguggenberger.
\newblock {\em Multiplication of Distributions and Applications to Partial
  Differential Equations}.
\newblock Pitman Research Notes in Mathematics 259. Longman Scientific {\&}
  Technical, 1992.

\end{thebibliography}
\newcommand{\SortNoop}[1]{}

\end{document}